\documentclass[11pt]{article}

\usepackage[applemac]{inputenc} 

\usepackage{amsthm,amssymb,amsbsy,amsmath,amsfonts,amssymb,amscd,mathrsfs}

\long\def\symbolfootnote[#1]#2{\begingroup%
\def\thefootnote{\fnsymbol{footnote}}\footnote[#1]{#2}\endgroup} 

\usepackage{graphicx}
\DeclareGraphicsExtensions{.jpg}
\usepackage{color}
\usepackage{fullpage}
\usepackage[hypertex]{hyperref}
\newtheorem{theorem}{Theorem}
\newtheorem{corollary}[theorem]{Corollary}
\newtheorem{lemma}[theorem]{Lemma}
\newtheorem{proposition}[theorem]{Proposition}
\theoremstyle{definition} 
\newtheorem{definition}[theorem]{Definition}
\newtheorem{example}[theorem]{Example}
\theoremstyle{remark}
\newtheorem{remark}[theorem]{Remark}

\newcommand{\bt}{\begin{theorem}}
\newcommand{\et}{\end{theorem}}
\newcommand{\bl}{\begin{lemma}}
\newcommand{\el}{\end{lemma}}
\newcommand{\bp}{\begin{proposition}}
\newcommand{\ep}{\end{proposition}}
\newcommand{\bc}{\begin{corollary}}
\newcommand{\ec}{\end{corollary}}
\newcommand{\bdeff}{\begin{definition}}
\newcommand{\edeff}{\end{definition}}
\newcommand{\brem}{\begin{remark}}
\newcommand{\erem}{\end{remark}}
\newcommand{\bex}{\begin{example}}
\newcommand{\eex}{\end{example}}

\newcommand{\con}{{\mathcal C}}
\newcommand{\mc}[1]{\mathcal {#1}}
\newcommand{\mb}[1]{\mathbb {#1}}

\newcommand{\bi}{\begin{itemize}}
\newcommand{\iii}{\item}
\newcommand{\ei}{\end{itemize}}
\newcommand{\bd}{\begin{description}}
\newcommand{\ed}{\end{description}}
\newcommand{\all}{\forall \,}
\newcommand{\bqn}{\begin{equation}}    
\newcommand{\eqn}{\end{equation}}       
\newcommand{\bqna}{\begin{eqnarray}}
\newcommand{\eqna}{\end{eqnarray}}
\newcommand{\eqnan}{\nonumber\end{eqnarray}}

\newcommand{\nn}{\nonumber}
\newcommand{\ba}[1]{\begin{array}{#1}}
\newcommand{\ea}{\end{array}}

\newcommand{\R}{\mathbb{R}}
\newcommand{\N}{\mathbb{N}}
\newcommand{\C}{\mathbb{C}}

\newcommand{\virg}[1]{``#1''}

\newcommand{\de}[2]{\partial_{#2} #1}

\newcommand{\lam}{\lambda}
\newcommand{\g}{\gamma}
\newcommand{\al}{\alpha}
\newcommand{\eps}{\varepsilon}

\newcommand{\tx}[1]{\mathrm{#1}}
\newcommand{\til}[1]{\widetilde{#1}}
\newcommand{\la}{\left\langle}
\newcommand{\ra}{\right\rangle}

\newcommand{\distr}{\Delta}
\newcommand{\metr}{{\bf g}}
\newcommand{\hh}{{\cal H}^{Q}}

\newcommand{\sss}{{\cal S}^{Q}}
\newcommand{\sub}{\mathbf{S}}
\newcommand{\ga}{\gamma}
\newcommand{\Pg}[1]{\left\{ #1 \right\}}

\newcommand{\hp}{hypothesis}
\newcommand{\EXP}{\textsf{Exp}}
\newcommand{\Exp}{\textsf{Exp}}

\newcommand{\Lam}{\Lambda}

\begin{document}
\begin{center} \noindent
{\LARGE{\sl{\bf On the Hausdorff volume in sub-Riemannian geometry}}}
\vskip 0.6 cm
Andrei Agrachev\\ 
{\footnotesize SISSA, Trieste, Italy and MIAN, Moscow, Russia - {\tt agrachev@sissa.it}}\\
\vskip 0.3 cm
Davide Barilari\\ 
{\footnotesize SISSA, Trieste, Italy - {\tt barilari@sissa.it}}\\
\vskip 0.3cm
Ugo Boscain\symbolfootnote[0]{This research has been supported by the European Research Council, ERC StG 2009 \virg{GeCoMethods}, contract number 239748, by the ANR Project GCM, program \virg{Blanche}, project number NT09-504490 and by the DIGITEO project CONGEO.}\\
{\footnotesize CNRS, CMAP Ecole Polytechnique, Paris, France - {\tt boscain@cmap.polytechnique.fr}}
\vskip 0.3cm
\end{center}
\vskip 0.1 cm

\begin{center}
\today
\end{center}
\vskip 0.2 cm
\begin{abstract} 
For a regular sub-Riemannian manifold we study the Radon-Nikodym derivative of the spherical Hausdorff measure with respect to a smooth volume. We prove that this is the volume of the unit ball in the nilpotent approximation and it is always a continuous function. We then prove that up to dimension 4 it is smooth, while starting from dimension 5, in corank 1 case, it is $\con^{3}$ (and $\con^{4}$ on every smooth curve) but in general not $\con^{5}$. These results answer to a question addressed by Montgomery about the relation between two intrinsic volumes that can be defined in a sub-Riemannian manifold, namely the Popp and the Hausdorff volume. If the nilpotent approximation depends on the point (that may happen starting from dimension 5), then they are not proportional, in general. 
\end{abstract}

\section{Introduction}

In this paper, by  a sub-Riemannian manifold we mean a triple $\sub=(M,\distr,\metr)$, where $M$ is a connected orientable smooth manifold of dimension $n$, $\distr$ is a smooth vector distribution of constant rank $k<n$, satisfying the H\"ormander condition and $\metr$ is an Euclidean structure on $\distr$.

A sub-Riemannian manifold has a natural structure of metric space, where the distance is the so called Carnot-Caratheodory distance
\begin{eqnarray}
d(q_0,q_1)=
\inf\{\int_0^T\sqrt{\metr_{\g(t)}(\dot\gamma(t),\dot\gamma(t))}~dt~|~ \gamma:[0,T]\to M \mbox{ is a Lipschitz curve},\cr
\gamma(0)=q_0,\gamma(T)=q_1, ~~\dot \gamma(t)\in\distr_{\gamma(t)}\mbox{ a.e. in $[0,T]$} \}.\nonumber
\end{eqnarray}
As a consequence of the H\"ormander condition this distance is always finite and continuous, and induces on $M$ the original topology (see \cite{chow,rashevsky}).

Since $(M,d)$  is a metric space, for every $\alpha>0$ one can define the $\alpha$-dimensional Hausdorff measure on $M$, and compute the Hausdorff dimension of $M$.

Define $\distr^{1}:=\distr, \distr^{i+1}:=\distr^{i}+[\distr^{i},\distr]$, for every $i=1,2,\ldots$ Under the hypothesis that the sub-Riemannian manifold is regular, i.e. 
if the dimension of $\distr^{i},\  i=1,\ldots,m$ do not depend on the point, the H\"ormander condition guarantees that there exists (a mimimal) $m\in \N$, called \emph{step} of the structure, such that $\distr_{q}^{m}=T_{q}M$, for all $q\in M$. The sequence
$$\mc{G}(\sub):=
(\underset{\begin{smallmatrix} \shortparallel \\ k
\end{smallmatrix}
}{\text{dim}\,\distr},\text{dim}\,\distr^{2},\ldots,\underset{\begin{smallmatrix} \shortparallel \\ n
\end{smallmatrix}
}{\text{dim}\, \distr^{m}})$$
is called \emph{growth vector} of the structure.

In this case, the graded vector space associated to the filtration $\distr_{q}\subset\distr_{q}^{2}\subset \ldots \subset \distr_{q}^{m}=T_{q}M$, 
$$\text{gr}_{q}(\distr)=\bigoplus_{i=1}^{m} \distr^{i}_{q}/\distr^{i-1}_{q}, \qquad\tx{where} \ \  \distr_{q}^{0}=0.$$
is well defined.
Moreover, it is well known that the Hausdorff dimension of $M$ is given by the formula (see \cite{mitchell})

$$Q=\sum_{i=1}^m i k_i, \qquad k_i:=\tx{dim}\, \distr_{q}^i/\distr_{q}^{i-1}.$$
In particular the Hausdorff dimension is always bigger than the topological dimension of $M$.

 Moreover, the $Q$-dimensional Hausdorff measure (denoted by $\hh$ in the following) behaves like a volume.
More precisely, in \cite{mitchell} Mitchell proved that if $\mu$ is a 
smooth volume\footnote{In the following by a smooth volume on $M$ we mean a measure $\mu$ associated to a smooth non-vanishing $n$-form $\omega_{\mu} \in \Lambda^{n}M$, i.e. for every measurable subset $A\subset M$ we set
$\mu(A)=\int_{A}\omega_{\mu}.$}
 on $M$ , then $d\mu=f_{\mu{\cal H}} d\hh$, where $f_{\mu{\cal H}}$ is a positive measurable function that is locally bounded and locally bounded away from zero, that is the Radon-Nikodym derivative of $\mu$ with respect to $\mc{H}^{Q}$.
 According to Mitchell terminology, this means that the two measures are \emph{commensurable} one with respect to the other.

Notice that the Hausdorff measure on sub-Riemannian manifolds has been intensively studied, see for instance \cite{gromov, mitchell}. A deep study of the Hausdorff measure for hypersurfaces in sub-Riemannian geometry, in particular in the context of Carnot groups, can be found in \cite{ambrosioserra,balogh,baloghtyson,lanconelli,capogna,franchiserra,magnani,serramonti} and references therein. Hausdorff measures for curves in sub-Riemannian manifolds were also studied in the problem of motion planning and complexity, see \cite{jp3,jp1,jp2,jean}.

Let us recall that there are  two common non-equivalent definitions of Hausdorff measure. The standard Hausdorff measure, 
where arbitrary coverings can be used, and the spherical Hausdorff measure,
where only ball-coverings appear (see Definition \ref{def:hauss}). 
However it is well known that, if $\sss$ denotes the spherical $Q$-Hausdorff measure, then $\hh$ is commensurable with $\sss$.\footnote{Indeed they are 
absolutely continuous one with respect to the other. In particular, for every $\alpha>0$, we have
$
2^{-\alpha}\mc{S}^{\alpha}\leq \mc{H}^{\alpha}\leq  \mc{S}^{\alpha}$ (see for instance \cite{federer}).
}
As a consequence, $\sss$ is commensurable with $\mu$, i.e.$$
d\mu=f_{\mu \mc{S}}d\sss,
$$
for a  positive measurable function $f_{\mu{\cal S}}$ that is locally bounded and locally bounded away from zero.
In this paper, we are interested to the properties of the 
function $f_{\mu{\cal S}}$. In particular, we would like to get informations about its regularity.

The reason why we study the spherical Hausdorff measure and not the standard Hausdorff measure  is that the first one appears to be more natural. Indeed, as explained later, $f_{\mu{\cal S}}$ is determined by the volume of the unit sub-Riemannian ball of the nilpotent approximation of the sub-Riemannian manifold, that can be explicitly described in a certain number of cases (see Theorem \ref{t-density} below).  
On the other hand nothing is known on how to compute $f_{\mu{\cal H}}$.
We conjecture that  $f_{\mu {\cal H}}$ is given by the $\mu$-volume of certain isodiametric sets, i.e. the maximum of the $\mu$-volume among all sets of diameter 1 in the nilpotent approximation (see \cite{leonardirigot,rigot} and reference therein for a discussion on isodiametric sets). This quantity is not very natural in sub-Riemannian geometry and is extremely difficult to compute. 

Our interests in studying $f_{\mu{\cal S}}$ comes from the following question:
\bd
\iii[Q1] How can we define an intrinsic volume in a sub-Riemannian manifold?
\ed
Here by intrinsic we mean a volume  which depends neither on the choice of the coordinate system, nor on the choice of the orthonormal frame, but only on the sub-Riemannian structure.

This question was first pointed out by Brockett, see \cite{brockett}, and by Montgomery in his book \cite{montgomerybook}. Having a volume that depends only on the geometric structure is interesting by itself, however, it is also necessary to  define intrinsically 
a Laplacian in a sub-Riemannian manifold. We recall that the Laplacian is defined as the divergence of the gradient and the definition of the divergence needs a volume since it measures how much the flow of a  vector field increases or decreases the volume.

Before talking about the question {\bf Q1} in sub-Riemannian geometry, let us briefly discuss it in the Riemannian case. In a 
$n$-dimensional Riemannian manifold there are three common ways of defining an invariant volume. The first is defined through the Riemannian structure and it is the so called Riemannian volume, which in coordinates has the expression $\sqrt{g}\, dx^{1}\ldots dx^n$, where $g$ is the determinant of the metric. The second and the third ones are defined via the Riemannian distance and are the $n$-dimensional Hausdorff measure and the $n$-dimensional spherical Hausdorff measure. These three volumes are indeed proportional (the constant of proportionality depending on the normalization, see e.g. 
\cite{chaveliso,federer}).

For what concern sub-Riemannian geometry, a regular sub-Riemannian manifold is a metric space, hence it is possible to define the Hausdorff volume $\hh$ and the spherical Hausdorff volume $\sss$. Also, there is an equivalent of the Riemannian volume, the so called Popp's volume $\cal{P}$, introduced by Montgomery in his book \cite{montgomerybook} (see also \cite{hypoelliptic}). The Popp volume is a smooth volume and was used in \cite{hypoelliptic} to define intrinsically the Laplacian (indeed a sub-Laplacian) in sub-Riemannian geometry.

In his book, Montgomery proposed to study whether these invariant  volumes are proportional as it occours in Riemannian geometry. More precisely, he addressed the following question:


\bd

\iii[Q2] Is Popp's measure equal to a constant multiple (perhaps depending on the growth vector) of the Hausdorff measure?
\ed
Mongomery noted that the answer to this question is positive for  left-invariant sub-Riemannian structures on  Lie groups, since  the Hausdorff (both the standard and spherical one) and the Popp volumes are left-invariant and hence proportional to the left Haar measure. But this question is nontrivial when there is no group structure.

One of the main purpose of our analysis is to answer to question {\bf Q2} for the spherical Hausdorff measure, i.e. to the question if the function $f_{{\cal PS}}$ (defined by $d\mc{P}=f_{\mc{PS}}d\mc{S}^{Q}$) is constant or not. More precisely,  we get a positive answer  for regular sub-Riemannian manifolds of dimension $3$ and $4$, while a negative answer starting from dimension 5, in general. 


Once a negative answer to {\bf Q2} is given,  it is natural to ask 
\bd
\iii[Q3] What is the regularity of $f_{{\cal PS}}$?
\ed
This question is important since the definition of an intrinsic Laplacian via $\sss$ require  $f_{ {\cal PS}}$ to be at least $\con^1$. 

Notice that since the Popp measure is a smooth volume, then $f_{\mu\mc{S}}$ is $\con^{k}, k=0,1,\ldots,\infty$ if and only if $f_{\mc{PS}}$ is as well.

We prove that $f_{\mu\mc{S}}$ is a continuous function and that for $n \leq 4$ it is smooth. In dimension 5 it is $\con^3$ but not smooth, in general. Moreover, we prove that the same result holds in all corank 1 cases (see Section \ref{s:bd} for a precise definition).

Our main tool is the nilpotent approximation (or the symbol) of the sub-Riemannian structure.
Recall that, under the regularity hypothesis, the sub-Riemannian structure $\sub=(M,\distr,\metr)$ induces a structure of nilpotent Lie algebra on $\tx{gr}_{q}(\distr)$. The nilpotent approximation at $q$ is the nilpotent simply connected Lie group $\tx{Gr}_{q}(\distr)$ generated by this Lie algebra, endowed with a suitable left-invariant sub-Riemannian structure $\widehat{\sub}_{q}$ induced by $\sub$, as explained in Section \ref{s:nilpapp}.

Recall that there exists a canonical isomorphism of 1-dimensional vector spaces (see \cite{hypoelliptic} for details)
\bqn \label{eq:iso}\bigwedge^{n}(T^{*}_{q}M) \simeq \bigwedge^{n}(\tx{gr}_{q}(\distr)^{*}).\eqn
Given a smooth volume $\mu$ on $M$, we define
the induced volume $\widehat{\mu}_{q}$ on the nilpotent approximation at point $q$ as the left-invariant volume on $\tx{Gr}_{q}(\distr)$ canonically associated to $\omega_{\mu}(q)\in \wedge^{n}(T^{*}_{q}M)$ by the above isomorphism. 

The first result concerns an explicit formula for $f_{\mu \mc{S}}$.
\bt
\label{t-density} Let  $\sub=(M,\distr,\metr)$ be a regular sub-Riemannian manifold. Let $\mu$ a volume on $M$ and $\widehat{\mu}_q$ the induced volume on the nilpotent approximation at point $q\in M$. If $A$ is an open subset of $M$, then
\bqn \label{eq:formulabella}
\mu(A)=\frac{1}{2^Q}\int_A \widehat{\mu}_q(\widehat{B}_q)\, d\mathcal{S}^Q,
\eqn
where $\widehat{B}_q$ is the unit ball in the nilpotent approximation at point $q$, i.e. 
$$f_{\mu \mc{S}}(q)= \frac{1}{2^{Q}}\widehat{\mu}_q(\widehat{B}_q).$$
\et

Starting from this formula we prove our first result about regularity of the density:
\bc\label{c-cont}
Let  $\sub=(M,\distr,\metr)$ be a regular sub-Riemannian manifold and let $\mu$ be a smooth volume on $M$. Then  the density $f_{\mu \mc{S}}$ is a continuous function.
\ec
Theorem \ref{t-density}, specified for the Popp measure $\mc{P}$, permits to answer the Montgomery's question. 
Indeed, the measure $\widehat{\mc{P}}_{q}$ induced by $\mc{P}$ on the nilpotent approximation at point $q$ coincides with the Popp measure built on $\widehat{\sub}_{q}$, as a sub-Riemannian structure. In other words, if we denote $\mc{P}_{\widehat{q}}$ the Popp measure on $\widehat{\sub}_{q}$, we get 
\bqn \label{eq:popp}\widehat{\mc{P}}_{q}=\mc{P}_{\widehat{q}}.\eqn

Hence, if the nilpotent approximation does not depend on the point, then $f_{{\cal PS}}$ is constant. In other words we have the following corollary.
\bc
Let  $\sub=(M,\distr,\metr)$ be a regular sub-Riemannian manifold and $\widehat{\sub}_{q}$ its nilpotent approximation at point $q\in M$. If $\widehat{\sub}_{q_{1}}$ is isometric to $\widehat{\sub}_{q_{2}}$ for any $q_{1},q_{2}\in M$, then $f_{\mc{PS}}$ is constant. 
In particular this happens if the sub-Riemannian structure is free.
\ec \noindent
For the definition of free structure see \cite{montgomerybook}.

Notice that, in the Riemannian case, nilpotent approximations at different points are isometric, hence the Hausdorff measure is proportional to the Riemannian volume (see \cite{chaveliso,federer}).

When the nilpotent approximation contains parameters that are function of the point, then, in general,  $f_{{\cal PS}}$ is not constant.
We have analyzed in details all growth vectors in dimension less or equal than 5:
\bi
\iii[-] dimension 3: (2,3), 
\iii[-] dimension 4:  (2,3,4), (3,4), 
\iii[-] dimension 5: (2,3,5), (3,5), (4,5) and the non generic cases (2,3,4,5), (3,4,5).
\ei
In all cases the nilpotent approximation is unique, except for the (4,5) case.
 As a consequence, we get:
\bt
\label{t-nleq5}
Let  $\sub=(M,\distr,\metr)$ be a regular sub-Riemannian manifold of dimension $n\leq 5$. Let $\mu$ be a smooth volume on $M$ and $\mc{P}$ be the Popp measure. Then 
\bi
\iii[$(i)$]if $\mc{G}(\sub)\neq (4,5)$, then $f_{\mc{PS}}$ is constant. As a consequence $f_{\mu \mc{S}}$ is smooth.
\iii[$(ii)$]if $\mc{G}(\sub)= (4,5)$, then $f_{\mu \mc{S}}$ is $\con^{3}$ (and $\con^{4}$ on smooth curves) but not $\con^{5}$, in general.
\ei
\et
Actually the regularity result obtained in the (4,5) case holds for all corank 1 structures, as specified by the following theorem.
\bt
\label{t-c4c5}
Let  $\sub=(M,\distr,\metr)$ be a regular corank 1 sub-Riemannian manifold of dimension $n\geq 5$. Let $\mu$ be a smooth volume on $M$. Then $f_{\mu\mc{S}}$ is $\con^{3}$ (and $\con^{4}$ on smooth curves) but not $\con^{5}$, in general.
\et
\noindent
Recall that for a corank 1 structure one has $\mc{G}(\sub)=(n-1,n)$ (see also Section \ref{s:bd}).

Notice that Theorem \ref{t-c4c5} apply in particular for the Popp measure.
The loss of regularity of $f_{\mu{\cal S}}$ is due to the presence of what are called \emph{resonance points}.  More precisely, the parameters appearing in the nilpotent approximation are the eigenvalues of a certain skew-symmetric matrix which depends on the point. Resonances are the points in which these eigenvalues are multiple.

To prove Theorem \ref{t-c4c5}, we have computed explicitly the optimal synthesis (i.e. all curves that minimize distance starting from one point) of the nilpotent approximation and, as a consequence,  the volume of nilpotent balls $\widehat{B}_q$.

Another byproduct of our analysis is 
\bp \label{p}
Under the hypothesis of Theorem \ref{t-c4c5}, if there are no resonance points then $f_{\mu \mc{S}}$ is smooth.
\ep

The structure of the paper is the following. In Section \ref{s:bd} we recall basic facts about sub-Riemannian geometry and about Hausdorff measures. In Section \ref{s:nilpapp} we provide normal forms for nilpotent structures in dimension less or equal than 5. In Section \ref{s:proof1} we prove Theorem \ref{t-density} and its corollaries, while in Section \ref{s:diff} we study the differentiability of the density for the corank 1 case. In the last Section we prove Theorem \ref{t-nleq5}.

\section{Basic Definitions}\label{s:bd}
\subsection{Sub-Riemannian manifolds}
We start recalling the definition of sub-Riemannian manifold.
\bdeff
A \emph{sub-Riemannian manifold} is a triple $\sub=(M,\distr,{\mathbf g})$, 
where
\bi
\iii[$(i)$] $M$ is a connected orientable smooth manifold of dimension $n\geq 3$;
\iii[$(ii)$] $\distr$ is a smooth distribution of constant rank $k< n$ satisfying the \emph{H\"ormander condition}, i.e. a smooth map that associates to $q\in M$  a $k$-dimensional subspace $\distr_{q}$ of $T_qM$ and we have
\bqn \label{Hor}
\text{span}\{[X_1,[\ldots[X_{j-1},X_j]]](q)~|~X_i\in\overline{\distr},\, j\in \N\}=T_qM, \quad \all q\in M,
\eqn
where $\overline{\distr}$ denotes the set of \emph{horizontal smooth vector fields} on $M$, i.e. $$\overline{\distr}=\Pg{X\in\mathrm{Vec}(M)\ |\ X(q)\in\distr_{q}~\ \forall~q\in M}.$$
\iii[$(iii)$] $\mathbf{g}_q$ is a Riemannian metric on $\distr_{q}$ which is smooth 
as function of $q$. We denote  the norm of a vector $v\in \distr_{q}$ with $|v|$, i.e.  $|v|=\sqrt{\metr_{q}(v,v)}.$
\ei
\edeff

A Lipschitz continuous curve $\ga:[0,T]\to M$ is said to be \emph{horizontal} (or \emph{admissible}) if 
$$\dot\ga(t)\in\distr_{\ga(t)}\qquad \text{ for a.e. } t\in[0,T].$$

Given an horizontal curve $\ga:[0,T]\to M$, the {\it length of $\ga$} is
\bqn
\label{e-lunghezza}
l(\ga)=\int_0^T |\dot{\g}(t)|~dt.
\eqn
The {\it distance} induced by the sub-Riemannian structure on $M$ is the 
function
\bqn
\label{e-dipoi}
d(q_0,q_1)=\inf \{l(\ga)\mid \ga(0)=q_0,\ga(T)=q_1, \ga\ \mathrm{horizontal}\}.
\eqn
The \hp\ of connectedness of $M$ and the H\"ormander condition guarantees the finiteness and the continuity of $d(\cdot,\cdot)$ with respect to the topology of $M$ (Chow-Rashevsky theorem, see, for instance, \cite{agrachevbook}). The function $d(\cdot,\cdot)$ is called the \emph{Carnot-Caratheodory distance} and gives to $M$ the structure of metric space (see \cite{bellaiche,gromov}).

\brem \label{r:action} It is a standard fact that $l(\g)$ is invariant under reparameterization of the curve $\g$.
Moreover, if an admissible curve $\ga$ minimizes the so-called {\it action functional}
$$ J(\g):=\frac{1}{2}\int_0^T |\dot{\g}(t)|^{2}dt. $$
with $T$ fixed (and fixed initial and final point), then $|\dot{\g}(t)|$ is constant
and $\ga$ is also a minimizer of $l(\cdot)$.
On the other side, a minimizer $\ga$ of $l(\cdot)$ such that  $|\dot{\g}(t)|$ is constant is a minimizer of $J(\cdot)$ with $T=l(\ga)/v$.
\erem

Locally, the pair $(\distr,{\mathbf g})$ can be given by assigning a set of $k$ smooth vector fields spanning $\distr$ and that are orthonormal for ${\mathbf g}$, i.e.  
\bqn
\label{trivializable}
\distr_{q}=\text{span}\{X_1(q),\dots,X_k(q)\}, \qquad \qquad \metr_q(X_i(q),X_j(q))=\delta_{ij}.
\eqn
In this case, the set $\Pg{X_1,\ldots,X_k}$ is called a \emph{local orthonormal frame} for the sub-Riemannian structure. 





\bdeff
Let $\distr$ be a distribution. Its \emph{flag} is the sequence of distributions $\distr^{1}\subset\distr^{2}\subset\ldots$ defined through the recursive formula
$$\distr^{1}:=\distr,~~~~~~~~\distr^{i+1}:=\distr^{i}+[\distr^{i},\distr].$$

A sub-Riemannian manifold is said to be \emph{regular} if for each $i=1,2,\ldots$ the dimension of $\distr^{i}_{q_{0}}$ does not depend on the point $q_0\in M$.
\edeff

\brem In this paper we always deal with regular sub-Riemannian manifolds. In this case H\"ormander condition can be rewritten as follows: 
$$\exists  \ \  \text{minimal}\  \  m\in \N \quad \text{ such that }\quad \distr^{m}_{q}=T_{q}M,\quad \all q\in M.$$
The sequence
$ \mathcal{G}(\sub):=(\text{dim}\,\distr,\text{dim}\,\distr^{2},\ldots,\text{dim}\, \distr^{m})$ is called \emph{growth vector}. 
  Under the regularity assumption $\mathcal{G}(\sub)$ does not depend on the point and $m$ is said the \emph{step} of the structure. The minimal growth is $(k,k+1,k+2,\ldots,n)$. When the growth is maximal the sub-Riemannian structure is called \emph{free} (see \cite{montgomerybook}).
\erem
A sub-Riemannian manifold is said to be \emph{corank 1} if its growth vector satisfies $\mathcal{G}(\sub)=(n-1,n)$.
A sub-Riemannian manifold $\sub$ of odd dimension is said to be \emph{contact} if 
$\distr= \ker \omega$, where $\omega \in \Lambda^{1}M$ and $d\omega|_{\distr} $ is non degenerate. 
A sub-Riemannian manifold $M$ of even dimension is said to be \emph{quasi-contact} if 
$\distr= \ker \omega$, where $\omega \in \Lambda^{1}M$ and satisfies $\tx{dim }\ker d\omega|_{\distr}=1$. 

Notice that contact and quasi-contact structures are regular and corank 1.

A sub-Riemannian manifold is said to be \emph{nilpotent} if there exists an orthonormal frame for the structure $\{X_{1},\ldots,X_{k}\}$ and $j\in \mb{N}$ such that
$[X_{i_{1}},[X_{i_{2}},\ldots,[X_{i_{j-1}},X_{i_{j}}]]]=0$ for every commutator of length $j$.

\subsection{Geodesics}\label{s:pmp}
 In this section we briefly recall some facts about sub-Riemannian geodesics. In particular, we define the sub-Riemannian Hamiltonian.

 \bdeff A \emph{geodesic} for a sub-Riemannian manifold $\sub=(M,\distr,{\mathbf g})$ is a curve $\ga:[0,T]\to M$ such that for every sufficiently small interval $[t_1,t_2]\subset [0,T]$, the restriction $\ga_{|_{[t_1,t_2]}}$ is a minimizer of $J(\cdot)$.
A geodesic for which $\mathbf{g}_{\ga(t)}(\dot \ga(t),\dot \ga(t))$  is (constantly) equal to one is said to be parameterized by arclength.\edeff

 Let us consider the cotangent bundle $T^{*}M$ with the canonical projection $\pi:T^{*}M\to M$, and denote the standard pairing between vectors and covectors with $\la\cdot,\cdot\ra$. 
The Liouville 1-form $s\in \Lambda^{1}(T^{*}M)$ is defined as follows: $s_{\lam}=\lam \circ \pi_{*}$, for every $\lam \in T^{*}M$. The canonical symplectic structure on $T^{*}M$ is defined by the closed 2-form $\sigma=ds$. In canonical coordinates $(\xi,x)$
 $$s= \sum_{i=1}^{n} \xi_{i} dx_{i}, \qquad \sigma = \sum_{i=1}^{n} d\xi_{i} \wedge dx_{i}.$$
We denote the Hamiltonian vector field associated to a function $h\in C^{\infty}(T^{*}M)$ with $\vec{h}$. Namely we have $dh=\sigma(\cdot,\vec{h})$ and in coordinates we have 
$$\vec{h}=\sum_i \frac{\partial h}{\partial \xi_i}\frac{\partial}{\partial x_i}-\frac{\partial h}{\partial x_i}\frac{\partial}{\partial \xi_i}$$
 
 The sub-Riemannian structure defines an Euclidean norm $|\cdot|$ on the distribution $\Delta_{q}\subset T_{q}M$. As a matter of fact this induces a dual norm  $$\|\lam\|= \max_{v\in \Delta_{q} \atop |v|=1} \la \lam,v\ra, \qquad \lam \in T_{q}^{*}M,$$  
 which is well defined on $\Delta_{q}^{*}\simeq T^{*}_{q}M/ \Delta_{q}^{\perp}$, 
 where $\Delta_{q}^{\perp}=\{\lam \in T_{q}^{*}M | \la\lam,v\ra=0, \all v\in \distr_{q}\}$ is the annichilator of the distribution.
 
 The \emph{sub-Riemannian Hamiltonian} is the smooth function on $T^{*}M$, which is quadratic on fibers, defined by
 $$H(\lam)=\frac{1}{2} \|\lam\|^{2}, \qquad \lam \in T_{q}^{*}M.$$
 If $\{X_{1},\ldots,X_{k}\}$ is a local orthonormal frame for the sub-Riemannian structure it is easy to see that
 $$H(\lam)=\frac{1}{2}\sum_{i=1}^{k} \la \lam ,X_{i}(q)\ra^{2}, \qquad \lam \in T_{q}^{*}M, \quad q=\pi(\lam).$$
 
Let $\sub=(M,\distr,{\mathbf g})$ be a sub-Riemannian manifold and fix $q_0\in M$. We define the \emph{endpoint map} (at time 1) as $$F:\mathcal{U}\to M, \quad F(\g)=\g(1),$$ where  $\mathcal{U}$ denotes  the set of admissible trajectories starting from $q_0$ and defined in $[0,1]$. If we fix a point $q_1\in M$, the problem of finding shortest paths from $q_0$ to $q_1$ is equivalent to the following one
\bqn \label{eq:min}
\min_{F^{-1}(q_1)}J(\g),
\eqn
 where $J$ is the action functional (see Remark \ref{r:action}). 
Then Lagrange multipliers rule implies that any $\g\in \mathcal{U}$ solution of \eqref{eq:min} satisfies one of the following equations
\begin{gather} 
\lam_{1} D_{\g}F=d_{\g}J, \label{eq:la1}\\
\lam_{1} D_{\g}F=0,\label{eq:la2}
\end{gather}
for some nonzero covector $\lam_{1} \in T^*_{\g(1)}M$ associated to $\g$. 
The following characterization is a corollary of Pontryagin Maximum Principle (PMP for short, see for instance \cite{agrachevbook,librougo,jurdjevicbook,pontrybook}): 
\bt \label{t-pmp} Let $\g$ be a minimizer. A nonzero covector $\lam_{1}$ satisfies \eqref{eq:la1} or \eqref{eq:la2}  if and only if there exists a Lipschitz curve $\lam(t)\in T^{*}_{\g(t)}M$, $t\in [0,1]$, such that $\lam(1)=\lam_{1}$ and 
\bi
\iii[-] if \eqref{eq:la1} holds, then $\lam(t)$ 
is a solution of $\dot{\lam}(t)=\overrightarrow{H}(\lam(t))$ for a.e. $t\in [0,1]$,
\iii[-] if \eqref{eq:la2} holds, then $\lam(t)$ satisfies $\sigma(\dot{\lam}(t), T_{\lam(t)}\distr^{\perp})=0$ for a.e. $t\in [0,1]$.
\ei
The curve $\lam(t)$ is said to be an extremal associated to $\gamma(t)$. In the first case $\lam(t)$ is called a \emph{normal extremal} while in the second one an \emph{abnormal extremal}. 
 \et
\brem \label{r-1/2} 
It is possible to give a unified characterization of normal and abnormal extremals in terms of the symplectic form. Indeed the Hamiltonian $H$ is always constant on extremals, hence $\lam(t)\subset H^{-1}(c)$ for some $c\geq0$.  Theorem \ref{t-pmp} can be rephrased as follows: any extremal $\lam(t)$ such that $H(\lam(t))=c$ is a reparametrization  of a characteristic curve of the differential form $\sigma|_{H^{-1}(c)}$, where $c=0$ for abnormal extremals, and $c>0$ for normal ones.

Also notice that, if $\lam(t)$ is a normal extremal, then, for every $\alpha>0$, $\lam_{\al}(t):=\alpha\, \lam (\alpha t)$ is also a normal extremal. If the curve is parametrized in such a way that $H=\frac{1}{2}$ then we say that the extremal is arclength parameterized. Trajectories parametrized by arclength corresponds to initial covectors $\lam_0$
belonging to the hypercylinder $\Lambda_{q_0}:=T^{*}_{q_{0}}M \cap H^{-1}(\frac{1}{2})\simeq S^{k-1}\times\R^{n-k}$ in  $T_{q_0}^\ast M$.

\erem

\brem \label{r:abnormal} 
From Theorem \ref{t-pmp} it follows that $\lam(t)=e^{t\vec{H}}(\lam_{0})$ is the normal extremal with initial covector $\lam_{0}\in \Lambda_{q_{0}}$.
If $\pi:T^{*}M \to M$ denotes the canonical projection, then it is well known that $\g(t)=\pi(\lam(t))$ is a  geodesic (starting from $q_{0}$). On the other hand, in every 2-step sub-Riemannian manifold all geodesics are projection of normal extremals, since there is no strict abnormal minimizer (see Goh conditions, \cite{agrachevbook}). 
\erem

The following proposition resumes some basic properties of small sub-Riemannian balls
\bp\label{p:2eps}
Let $\sub$ be a sub-Riemannian manifold and $B_{q_{0}}(\eps)$ the sub-Riemannian ball of radius $\eps$ at fixed point $q_{0}\in M$. For $\eps>0$ small enough we have:
\bi
\iii[(i)]  $\all q \in B_{q_{0}}(\eps)$ there exists a minimizer that join $q$ and $q_{0}$,
\iii[(ii)] $\tx{diam}(B_{q_{0}}(\eps))=2\eps$.
\ei
\ep
Claim $(i)$ is a consequence of Filippov theorem (see \cite{agrachevbook,bressanbook}). To prove $(ii)$ it is sufficient to show that, for $\eps$ small enough, there exists two points in $q_{1},q_{2}\in \partial B_{q_{0}}(\eps)$ such that $d(q_{1},q_{2})=2\eps$. 

To this purpose, consider the projection $\g(t)=\pi(\lam(t))$ of a normal extremal starting from $\g(0)=q_{0}$, and defined in a small neighborhood of zero $t\in]-\delta,\delta[$ . Using arguments of Chapter 17 of \cite{agrachevbook} one can prove that $\g(t)$ is globally minimizer. Hence if we consider $0<\eps<\delta$ we have that $q_{1}=\g(-\eps)$ and $q_{2}=\g(\eps)$ satisfy the property required, which proves claim $(ii)$.

\bdeff
Fix $q_0\in M$. We define the \emph{Exponential map} starting from $q_{0}$ as
$$\EXP_{q_{0}}: T^{*}_{q_{0}}M \to M, \qquad \EXP_{q_{0}}(\lam_{0})= \pi(e^{\vec{H}}(\lam_{0})).$$
Using the homogeneity property $H(c\lambda)=c^{2}H(\lambda), \ \all c>0$, we have that
$$e^{\vec{H}}(s\lam)=e^{s\vec{H}}(\lam), \qquad \all s>0.$$
In other words we can recover the geodesic on the manifold with initial covector $\lam_{0}$ as the image under $\EXP_{q_{0}}$ of the ray $\{t\lam_{0},t\in[0,1]\}\subset T^{*}_{q_{0}}M$ that join the origin to $\lam_{0}$.
$$\EXP_{q_{0}}(t\lam_{0})=\pi(e^{\vec{H}}(t\lam_{0}))= \pi(e^{t\vec{H}}(\lam_{0}))=\pi(\lam(t))=\g(t).$$
\edeff
%
Next, we recall the definition of cut and conjugate time. 
\bdeff \label{def:cutconj} Let $\sub=(M,\distr,\metr)$ be a sub-Riemannian manifold. 
Let $q_{0}\in M$ and $\lam_0\in\Lam_{q_0}$. Assume that the geodesic $\g(t)=\textsf{Exp}_{q_{0}}(t\lam_{0})$ for $t>0$, is not abnormal.
\bi
\iii[$(i)$] The \emph{first conjugate time} is 
$t(\lam_{0})=\min\{t>0,\  t\lam_{0}\  \text{is a critical point of}\  \EXP_{q_{0}}\}$.

\iii[$(ii)$] The \emph{cut time} is $t_{c}(\lam_{0})=\min\{t>0,\ \exists\, \lam_{1}\in \Lambda_{q_{0}},\lam_{1}\neq\lam_{0} \ \text{s.t.}\ \EXP_{q_{0}}(t_{c}(\lam_{0})\lam_{0})=\EXP_{q_{0}}(t_{c}(\lam_{0})\lam_{1}) \}.$
\ei
\edeff
It is well known that if a geodesic is not abnormal then it loses optimality either at the cut or at the conjugate locus (see for instance \cite{agrexp}).

\subsection{Hausdorff measures}\label{sec:misure}
In this section we recall definitions of Hausdorff measure and spherical Hausdorff measure. We start with the definition of smooth volume.
\bdeff  \label{17}Let $M$ be a $n$-dimensional smooth manifold, which is connected and orientable. By a \emph{smooth volume} on $M$ we mean a measure $\mu$ on $M$ associated to a smooth non-vanishing $n$-form $\omega_{\mu} \in \Lambda^{n}M$, i.e. for every subset $A\subset M$ we set
$$\mu(A)=\int_{A}\omega_{\mu}.$$
\edeff
The Popp volume $\mc{P}$, which is a smooth volume in the sense of Definition \ref{17}, is the volume associated to a $n$-form $\omega_{\mc{P}}$ that can be intrinsically defined via the sub-Riemannian structure (see \cite{hypoelliptic,montgomerybook}).

Let $(M,d)$ be a metric space and
denote with $\mc{B}$ the set of balls in $M$.
\bdeff \label{def:hauss} Let $A$ be a subset of $M$ and $\alpha>0$. 

The \emph{$\alpha$-dimensional Hausdorff measure} of $A$ is
$$\mc{H}^{\alpha}(A):= \lim_{\delta\to 0} \mc{H}^{\alpha}_{\delta}(A),$$
where
$$\mc{H}^{\alpha}_{\delta}(A):=\inf\{ \sum_{i=1}^{\infty} \text{diam}(A_{i})^{\alpha}, A\subset \bigcup_{i=1}^{\infty} A_{i}, \text{diam}(A_{i})<\delta\}.$$

The \emph{$\alpha$-dimensional spherical Hausdorff measure} of $A$ is
$$\mc{S}^{\alpha}(A):= \lim_{\delta\to 0} \mc{S}^{\alpha}_{\delta}(A),$$ where $$\mc{S}^{\alpha}_{\delta}(A):=\inf\{ \sum_{i=1}^{\infty} \text{diam}(B_{i})^{\alpha}, A\subset \bigcup_{i=1}^{\infty} B_{i}, B_{i}\in \mc{B}, \text{diam}(B_{i})<\delta\}.$$

These two measures are commensurable since it holds (see \cite{federer})
\bqn \label{eq:hsh}
{2^{-\alpha}}\mc{S}^{\alpha}(A)\leq \mc{H}^{\alpha}(A)\leq  \mc{S}^{\alpha}(A), \qquad \all A \subset M.
\eqn

The \emph{Hausdorff dimension} of $A$ is defined as 
\bqn\label{eq:hdim}
\inf \{\alpha>0, \mc{H}^{\alpha}(A)=0\}=\sup\{\alpha>0, \mc{H}^{\alpha}(A)=+\infty\}.
\eqn
\edeff
Formula \eqref{eq:hsh} guarantees that Hausdorff dimension of $A$ does not change if we replace $\mc{H}^{\alpha}$ with $\mc{S}^{\alpha}$ in formula \eqref{eq:hdim}.

It is a standard fact that the Hausdorff dimension of a Riemannian manifold, considered as a metric space, coincides with its topological dimension. 
On the other side, we have the following

\bt \label{t:mitchell} Let $(M,\distr,\metr)$ be a regular sub-Riemannian manifold. Its Hausdorff dimension as a metric space is
 $$Q=\sum_{i=1}^m i k_i, \qquad k_i:=\tx{dim}\, \distr^{i}-\tx{dim}\, \distr^{i-1}.$$
Moreover $\mc{S}^{Q}$ is commensurable to a smooth volume $\mu$ on $M$, i.e. for every compact $K\subset M$ there exists $\alpha_{1},\alpha_{2}>0$ such that
\bqn \label{eq:comp} \alpha_{1} \mathcal{S}^Q\leq \mu\leq \alpha_{2} \mathcal{S}^Q.\eqn \et
This theorem was proved by Mitchell in \cite{mitchell}. In its original version it was stated  for the Lebesgue measure and the standard Hausdorff measure.

\section{The nilpotent approximation}\label{s:nilpapp}
In this section we briefly recall the concept of nilpotent approximation. For details see \cite{agrachevlocal,bellaiche}.
\subsection{Privileged coordinates}
Let $\sub=(M,\distr,\metr)$ be a sub-Riemannian manifold and $(X_{1},\ldots,X_{k})$ an orthonormal frame. Fix a point $q\in M$ and consider the flag of the distribution 
$\distr^{1}_{q}\subset\distr^{2}_{q}\subset\ldots\subset\distr^{m}_{q}$. Recall that $k_i=\tx{dim}\, \distr^{i}_{q}-\tx{dim}\, \distr^{i-1}_{q}$ for $i=1,\ldots,m$, and that $k_{1}+\ldots+k_{m}=n$.

Let $O_{q}$ be an open neighborhood of the point $q\in M$. We say that a system of coordinates $\psi: O_{q}\to \R^{n}$ is \emph{linearly adapted} to the flag if, in these coordinates, we have $\psi(q)=0$ and
$$\psi_{*}(\distr^{i}_{q})=\R^{k_{1}}\oplus \ldots \oplus \R^{k_{i}}, \qquad \all i=1,\ldots,m.$$

Consider now the splitting
$\R^{n}=\R^{k_{1}}\oplus \ldots \oplus \R^{k_{m}}$
and denote its elements $x=(x_{1},\ldots,x_{m})$ where $x_{i}=(x_{i}^{1},\ldots,x_{i}^{k_{i}})\in \R^{k_{i}}$. 
The space of all differential operators in $\R^{n}$ with smooth coefficients forms an associative
algebra with composition of operators as multiplication. The differential operators with
polynomial coefficients form a subalgebra of this algebra with generators $1, x_{i}^{j} ,\frac{\partial}{\partial x_{i}^{j}},$ where
$i=1,\ldots,m;\  j = 1,\ldots,k_{i}$. We define weights of generators as
$$\nu(1)=0, \qquad \nu(x_{i}^{j})=i, \qquad \nu(\frac{\partial}{\partial x_{i}^{j}})=-i,$$ and the weight of monomials
$$\nu(y_{1}\cdots y_{\alpha}\frac{\partial^{\beta}}{\partial z_{1} \cdots \partial z_{\beta}})=\sum_{i=1}^{\alpha}\nu(y_{i})-\sum_{j=1}^{\beta}\nu(z_{j}).$$
Notice that a polynomial differential operator homogeneous with respect to $\nu$ (i.e. whose monomials are all of same weight) is homogeneous with respect to dilations $\delta_{t}:\R^{n}\to \R^{n}$ defined by
\bqn \label{dil} \delta_{t}(x_{1},\ldots,x_{m})=(tx_{1},t^{2}x_{2},\ldots, t^{m}x_{m}), \qquad t>0.
\eqn
In particular for a homogeneous vector field $X$ of weight $h$ it holds $\delta_{t*}X=t^{-h}X$.
A smooth vector field $X\in \tx{Vec}(\R^{n})$, as a first order differential operator, can be written as
$$X=\sum_{i,j} a_{i}^{j}(x) \frac{\partial}{\partial x_{i}^{j}}$$
and considering its Taylor expansion at the origin we can write the formal expansion
$$X\approx \sum_{h=-m}^{\infty} X^{(h)}$$
where $X^{(h)}$ is the homogeneous part of degree $h$ of $X$ (notice that every monomial of a first order differential operator has weight not smaller than $-m$).
Define the filtration of  $\tx{Vec}(\R^{n})$ $$\mc{D}^{(h)}=\{X\in \tx{Vec}(\R^{n}): X^{(i)}=0, \all i<h\}, \qquad \ell\in \mb{Z}.$$
\bdeff A system of coordinates $\psi: O_{q}\to \R^{n}$ defined near the point $q$ is said \emph{privileged} for a sub-Riemannian structure $\sub$ if these coordinates are linearly adapted to the flag and such that $\psi_{*}X_{i}\in \mc{D}^{(-1)}$ for every $i=1,\ldots,k$.
\edeff

\bt Privileged coordinates always exists. Moreover there exist $c_{1},c_{2}>0$ such that in these coordinates, for all $\eps >0$ small enough, we have
\bqn \label{eq:ballbox}
c_{1}\,\tx{Box}(\eps) \subset B(q,\eps) \subset c_{2}\,\tx{Box}(\eps) , 
\eqn
where $\tx{Box}(\eps)=\{x\in \R^{n}, |x_{i}|\leq \eps^{i}\}$.
\et

Existence of privileged coordinates is proved in \cite{agrachevlocal,agrachevlie,bellaiche,bianchinistefani}. In the regular case the construction of privileged coordinates was also done in the context of hypoelliptic operators (see \cite{rothstein}). The second statement is known as \emph{Ball-Box theorem} and a proof can be found in \cite{bellaiche}. Notice however that privileged coordinates are not unique.

\bdeff Let $\sub=(M,\distr,\metr)$ be a regular sub-Riemannian manifold and $(X_{1},\ldots,X_{k})$ a local orthonormal frame near a point $q$. Fixed a system of privileged coordinates, we define the \emph{nilpotent approximation of $\sub$ near $q$}, denoted by $\widehat{\sub}_{q}$, the sub-Riemannian structure on $\R^{n}$ having $(\widehat{X}_{1}, \ldots, \widehat{X}_{k})$ as an orthonormal frame, where $\widehat{X}_{i}:=(\psi_{*}X_{i})^{(-1)}$.
\edeff
\brem \label{palla} It is well known that under the regularity hypothesis, $\widehat{\sub}_{q}$ is naturally endowed with a Lie group structure whose Lie algebra is generated by left-invariant vector fields $\widehat{X}_{1}, \ldots, \widehat{X}_{k}$. Moreover the sub-Riemannian distance $\widehat{d}$ in $\widehat{\sub}_{q}$ is homogeneous with respect to dilations $\delta_{t}$, i.e. $\widehat{d}(\delta_{t}(x),\delta_{t}(y))=t \,\widehat{d}(x,y)$. In particular, if $\widehat{B}_{q}(r)$ denotes the ball of radius $r$ in $\widehat{\sub}_{q}$, this implies $\delta_{t}(\widehat{B}_{q}(1))=\widehat{B}_{q}(t)$.
\erem

\bt \label{t:napp}
The nilpotent approximation $\widehat{\sub}_{q}$ of a sub-Riemannian structure $\sub$ near a point $q$ is the metric tangent space to $M$ at point $q$ in the sense of Gromov, that means
\bqn\label{eq:gh}
\delta_{1/\eps} B(q,\eps) \longrightarrow \widehat{B}_{q},
\eqn
where $\widehat{B}_{q}$ denotes the sub-Riemannian unit ball of the nilpotent approximation $\widehat{\sub}_{q}$.
\et
\brem Convergence of sets in \eqref{eq:gh} is intended in the Gromov-Hausdorff topology \cite{bellaiche,gromovbook}. In the regular case this theorem was proved by Mitchell in \cite{mitchellphd}. A proof in the general case can be found in \cite{bellaiche}. 
\erem

\bdeff Let $\sub=(M,\distr,\metr)$ be a regular sub-Riemannian manifold and  $\widehat{\sub}_{q}$ its nilpotent approximation near $q$. If $\mu$ a smooth volume on $M$, associated to the smooth non-vanishing $n$-form $\omega_{\mu}$, we define the \emph{induced volume} $\widehat{\mu}_{q}$ at the point $q$ as the left-invariant volume on $\widehat{\sub}_{q}$ canonically associated with $\omega_{\mu}(q)\in \wedge^{n}(T^{*}_{q}M)$ (cf. isomorphism \eqref{eq:iso}).
\edeff
From Theorem \ref{t:napp} and the relation\footnote{Notice that this formula is meaningful in privileged coordinates near $q$.} $\mu(\delta_{\eps}A)=\eps^{Q}\widehat{\mu}_{q}(A)+o(\eps^{Q})$ when $\eps\to 0$, one gets
\bc\label{c:nd} Let $\mu$ be a smooth volume on $M$ and $\widehat{\mu}_{q}$ the induced volume on the nilpotent approximation at point $q$. Then, for $\eps \to 0$, we have 
$$\mu(B(q,\eps))=\eps^{Q}\widehat{\mu}_{q}(\widehat{B}_{q})+o(\eps^{Q}). $$
\ec

\subsection{Normal forms for nilpotent approximation in dimension $\leq 5$}
In this section we provide normal forms for the nilpotent approximation of regular sub-Riemannian structures in dimension less or equal than 5.
One can easily shows that in this case the only possibilities for growth vectors are:
\bi
\iii[-] dim$(M)=3$:\qquad $\mc{G}(\sub)=(2,3),$
\iii[-] dim$(M)=4$:\qquad$\mc{G}(\sub)=(2,3,4)$ or $\mc{G}(\sub)=(3,4),$
\iii[-] dim$(M)=5$:\qquad  $\mc{G}(\sub)\in \{(2,3,4,5), (2,3,5), (3,5), (3,4,5), (4,5)\}.$
\ei
We have the following.

\bt\label{t:nf}
Let $\sub=(M,\distr,\metr)$ be a regular sub-Riemannian manifold and $\widehat{\sub}_{q}$ its nilpotent approximation near $q$. Up to a change of coordinates and a rotation of the orthonormal frame we have the following expression for the orthonormal frame of $\widehat{\sub}_{q}$:
\begin{description}
\iii[Case $n=3$] 
\bi \iii $\mc{G}(\sub)=(2,3)$. (Heisenberg)
\begin{align*}
\widehat{X}_{1}&=\partial_{1},\\
\widehat{X}_{2}&=\partial_{2}+x_{1}\partial_{3}.
\end{align*}\ei
\iii[Case $n=4$]\bi \iii $\mc{G}(\sub)=(2,3,4)$. (Engel)
\begin{align*}
\widehat{X}_{1}&=\partial_{1},~~~~~~~~\\
\widehat{X}_{2}&=\partial_{2}+x_{1}\partial_{3}+x_{1}x_{2}\partial_{4}.
\end{align*}
\iii $\mc{G}(\sub)=(3,4)$. (Quasi-Heisenberg)
\begin{align*}
\widehat{X}_{1}&=\partial_{1},\\
\widehat{X}_{2}&=\partial_{2}+x_{1}\partial_{4},\\
\widehat{X}_{3}&=\partial_{3}.
\end{align*}
\ei
\iii[Case $n=5$]
\bi 
\iii $\mc{G}(\sub)=(2,3,5)$. (Cartan) 
\begin{align*}
\widehat{X}_{1}&=\partial_{1},\\
\widehat{X}_{2}&=\partial_{2}+x_{1}\partial_{3}+\frac{1}{2}x_{1}^{2}\partial_{4}+x_{1}x_{2}\partial_{5}.
\end{align*}
\iii $\mc{G}(\sub)=(2,3,4,5)$. (Goursat rank 2)
\begin{align*}
\widehat{X}_{1}&=\partial_{1},\\
\widehat{X}_{2}&=\partial_{2}+x_{1}\partial_{3}+\frac{1}{2}x_{1}^{2}\partial_{4}+\frac{1}{6}x_{1}^{3}\partial_{5}.
\end{align*}
\iii $\mc{G}(\sub)=(3,5).$ (Corank 2)
\begin{align*}
\widehat{X}_{1}&=\partial_{1}-\frac{1}{2}x_{2}\partial_{4},\\
\widehat{X}_{2}&=\partial_{2}+\frac{1}{2}x_{1}\partial_{4}-\frac{1}{2}x_{3}\partial_{5},\\
\widehat{X}_{3}&=\partial_{3}+\frac{1}{2}x_{2}\partial_{5}.
\end{align*}
\iii $\mc{G}(\sub)=(3,4,5).$ (Goursat rank 3)
\begin{align*}
\widehat{X}_{1}&=\partial_{1}-\frac{1}{2}x_{2}\partial_{4}-\frac{1}{3}x_{1}x_{2}\partial_{5},\\
\widehat{X}_{2}&=\partial_{2}+\frac{1}{2}x_{1}\partial_{4}+\frac{1}{3}x_{1}^{2}\partial_{5},\\
\widehat{X}_{3}&=\partial_{3}. 
\end{align*}
\iii $\mc{G}(\sub)=(4,5).$ (Bi-Heisenberg)
\begin{align}\label{eq:45}
\widehat{X}_{1}&=\partial_{1}-\frac{1}{2}x_{2}\partial_{5},\notag\\
\widehat{X}_{2}&=\partial_{2}+\frac{1}{2}x_{1}\partial_{5},\notag\\
\widehat{X}_{3}&=\partial_{3}-\frac{\alpha}{2}x_{4}\partial_{5}, \qquad \alpha\in \R,\\
\widehat{X}_{4}&=\partial_{4}+\frac{\alpha}{2}x_{3}\partial_{5}.\notag
\end{align}
\ei
\end{description}
\et
\begin{proof} It is sufficient to find, for every such a structure, a basis of the Lie algebra such that the structural constants\footnote{Let $X_{1},\ldots, X_{k}$  be a basis of a Lie algebra $\mathfrak{g}$. The coefficients $c_{ij}^{\ell}$ that satisfy $[X_{i},X_{j}]=\sum_{\ell}c_{ij}^{\ell}X_{\ell}$ are called structural constant of $\mathfrak{g}$.} are uniquely determined by the sub-Riemannian structure.
We give a sketch of the proof for the $(2,3,4,5)$ and $(3,4,5)$ and $(4,5)$ cases. The other cases can be treated in a similar way. 

$(i)$. Let $\widehat{\sub}=(G,\distr,\metr)$ be a nilpotent $(3,4,5)$ sub-Riemannian structure. Since we deal with a left-invariant sub-Riemannian structure, we can identify the distribution $\Delta$ with its value at the identity of the group $\distr_{id}$. Let $\{e_{1},e_{2},e_{3}\}$ be a basis for  $\distr_{id}$, as a vector subspace of the Lie algebra. By our assumption on the growth vector  we know that 
\bqn \label{eq:aaa} \tx{dim}\  \tx{span}\{[e_{1},e_{2}],[e_{1},e_{3}],[e_{2},e_{3}]\}/ \distr_{id} =1. \eqn
In other words, we can consider the 
skew-simmetric mapping
\bqn \label{eq:liebr}\Phi(\cdot,\cdot):=[\cdot,\cdot]/\distr_{id}:\distr_{id} \times \distr_{id} \to T_{id}G/\distr_{id},\eqn
and condition \eqref{eq:aaa} implies that there exists a one dimensional subspace in the kernel of this map. Let $\widehat{X}_{3}$ be a normalized vector in the  kernel and consider its orthogonal subspace $D\subset \distr_{id}$ with respect to the Euclidean product on $\distr_{id}$. Fix an arbitrary orthonormal basis $\{X_{1},X_{2}\}$ of $D$ and set $\widehat{X}_{4}:=[X_{1},X_{2}]$. It is easy to see that $\widehat{X}_{4}$ does not change if we rotate the base $\{X_{1},X_{2}\}$ and there exists a choice of this frame, denoted $\{\widehat{X}_{1},\widehat{X}_{2}\}$, such that $[\widehat{X}_{2},\widehat{X}_{4}]=0$. Then set $\widehat{X}_{5}:=[\widehat{X}_{1},\widehat{X}_{4}]$. Therefore we found a canonical basis for the Lie algebra that satisfies the following commutator relations:
\begin{align*}
[\widehat{X}_{1},\widehat{X}_{2}]&=\widehat{X}_{4}, \qquad  \qquad [\widehat{X}_{1},\widehat{X}_{4}]=\widehat{X}_{5}, 
\end{align*} and all other commutators vanish.
A standard application of the Campbell-Hausdorff formula gives the coordinate expression above.

\smallskip
$(ii)$. Let us assume now that $\widehat{\sub}$ is a nilpotent $(2,3,4,5)$ sub-Riemannian structure. 
As before we identify the distribution $\Delta$ with its value at the identity and consider any orthonormal basis $\{e_{1},e_{2}\}$ for the 2-dimensional subspace $\distr_{id}$. By our assumption on $\mc{G}(\sub)$
\begin{gather} 
\tx{dim}\  \tx{span}\{e_{1},e_{2},[e_{1},e_{2}]\} =3 \notag\\
\tx{dim}\  \tx{span}\{e_{1},e_{2},[e_{1},e_{2}],[e_{1},[e_{1},e_{2}]],[e_{2},[e_{1},e_{2}]]\} =4.  \label{xdccc}
\end{gather}
As in $(i)$, it is easy to see that there exists a choice of the orthonormal basis on $\distr_{id}$, which we denote $\{\widehat{X}_{1},\widehat{X}_{2}\}$, such that $[\widehat{X}_{2},[\widehat{X}_{1},\widehat{X}_{2}]]=0$. From this property and the Jacobi identity it follows $[\widehat{X}_{2},[\widehat{X}_{1},[\widehat{X}_{1},\widehat{X}_{2}]]]=0$. Then we set $\widehat{X}_{3}:=[\widehat{X}_{1},\widehat{X}_{2}]$, $\widehat{X}_{4}=[\widehat{X}_{1},[\widehat{X}_{1},\widehat{X}_{2}]]$ and $\widehat{X}_{5}:=[\widehat{X}_{1},[\widehat{X}_{1},[\widehat{X}_{1},\widehat{X}_{2}]]]$. It is easily seen that \eqref{xdccc} implies that these vectors are linearly independent and give a canonical basis for the Lie algebra, with  the only nontrivial commutator relations:
\begin{align*}
[\widehat{X}_{1},\widehat{X}_{2}]=\widehat{X}_{3}, \qquad [\widehat{X}_{1},\widehat{X}_{3}]=\widehat{X}_{4},  \qquad [\widehat{X}_{1},\widehat{X}_{4}]=\widehat{X}_{5}. 
\end{align*} 

$(iii)$. In the case $(4,5)$ since $\dim\, T_{id}G/\distr_{id}=1$, the 
 map \eqref{eq:liebr} is represented by a single $4 \times 4$ skew-simmetric matrix $L$. By skew-symmetricity its eigenvalues are purely imaginary $\pm i b_{1},\pm i b_{2}$, one of which is different from zero. Assuming $b_{1}\neq 0$ we have that $\alpha=b_{2}/b_{1}$. Notice that the structure is contact if and only if $\alpha \neq 0$ (see also Section \ref{s-norm} for more details on the normal form).
\end{proof}

\brem 
Notice that, in the statement of Theorem \ref{t:nf}, in all other cases the nilpotent approximation does not depend on any parameter, except for the $(4,5)$ case.
As a consequence, up to dimension 5, the sub-Riemannian structure induced on the tangent space, and hence the Popp measure $\mc{P}$, does not depend on the point, except for the $(4,5)$ case.

In the $(4,5)$ case we have the following expression for the Popp's measure $$\mc{P}=\frac{1}{\sqrt{b_{1}^{2}+b_{2}^{2}}}\,dx_{1}\wedge \ldots \wedge dx_{5},$$
where $b_{1},b_{2}$ are the eigenvalues of the skew-simmetric matrix that represent the Lie bracket map.
\erem

Since the normal forms in Theorem \ref{t:nf} do not depend on the point, except when $\mc{G}(\sub) \neq (4,5)$, we have the following corollary 

\bc\label{c:isometrie}
Let $\sub=(M,\distr,\metr)$ be a regular sub-Riemannian manifold such that $\tx{dim}(M)\leq 5$ and $\mc{G}(\sub) \neq (4,5)$. Then if $q_{1},q_{2}\in M$ we have that $\widehat{\sub}_{q_{1}}$ is isometric to $\widehat{\sub}_{q_{2}}$ as sub-Riemannian manifolds.
\ec

\section{Proof of Theorem \ref{t-density}: the density is the volume of nilpotent balls}\label{s:proof1}
In this section we prove Theorem \ref{t-density}, i.e. 
\bqn \label{eq}
f_{\mu \mc{S}}(q)= \frac{1}{2^{Q}}\widehat{\mu}_q(\widehat{B}_q).
\eqn

It is well known that, being $\mu$ absolutely continuous with respect to $\mc{S}^{Q}$ (see Theorem \ref{t:mitchell}), the Radon-Nikodym derivative of $\mu$ with respect to $\mc{S}^{Q}$, namely $f_{\mu \cal{S}}$, can be computed almost everywhere as
$$\lim_{r\to 0} \frac{\mu(B(q,r))}{\mc{S}^{Q}(B(q,r))}.$$
By Corollary \ref{c:nd} we get
\begin{align*}
\frac{\mu(B(q,r))}{\mc{S}^{Q}(B(q,r))}&= \frac{r^{Q}\widehat{\mu}_{q}(\widehat{B}_{q})+o(r^{Q})}{\mc{S}^{Q}(B(q,r))}
=  \frac{\widehat{\mu}_{q}(\widehat{B}_{q})}{2^{Q}}\frac{2^{Q}r^Q}{\mc{S}^{Q}(B(q,r))}+\frac{o(r^{Q})}{\mc{S}^{Q}(B(q,r))}. \end{align*}

Then we are left to prove the following


\bl Let $A$ be an open subset of $M$. For $\mathcal{S}^Q$-a.e. $q\in A$ we have 
\bqn \label{eq:123}
\lim_{r \rightarrow 0} \frac{\mathcal{S}^Q(A\cap B(q,r))}{(2r)^Q}=1.
\eqn
\el
\begin{proof} In the following proof we make use of Vitali covering lemma \footnote{\textbf{Theorem.}(Vitali covering lemma, \cite{federer,ambrosiometric})
Let $E$ be a metric space, $B\subset E$ and $\alpha>0$ such that
$\mathcal{H}^{\alpha}(B)<\infty$, and let $\mc{F}$ be a
fine covering of $B$. Then
there exist a countable disjoint subfamily $\{V_i\}\subset \mc{F}$
such that $$\mc{H}^{\alpha}(B\setminus \bigcup V_i)=0.$$
We recall that $\mc{F}$ is a \emph{fine covering} of
$B$ if for every $x\in B$ and $\eps>0$ there exists $V\in \mc{F}$
such that $x\in V$ and $\tx{diam}(V)<\eps$.}
and we always assume that balls of our covering are small enough to satisfy property $(ii)$ of Proposition \ref{p:2eps}.

We prove that the set where \eqref{eq:123} exists and is
different from 1 has $\mathcal{S}^{Q}$-null measure. 

\noindent
$(i).$ First we show
$$\mathcal{S}^{Q}(E_{\delta})=0,\qquad \forall\, 0<\delta\leq1,$$ where
$$E_\delta:=\{q\in A: \mathcal{S}^{Q}(A\cap B(q,r))<(1-\delta)(2r)^Q, \forall\,
0<r<\delta\}$$
Let $\{B_i\}$ a ball covering of $E_\delta$ with
$\tx{diam}(B_i)<\delta$ and such that
$$\sum_i \tx{diam}(B_i)^Q\leq \mathcal{S}^Q_{\delta}(E_{\delta})+\eps \leq \mathcal{S}^Q(E_{\delta})+\eps.$$
Then we have
\begin{align*}
\mathcal{S}^{Q}(E_{\delta})&\leq \mathcal{S}^{Q}(A \cap \bigcup B_i)
                          \\&\leq  \sum \mathcal{S}^{Q}(A \cap B_i) \\
                          &\leq  (1-\delta) \sum \tx{diam}(B_i)^Q 
                          \\&\leq (1-\delta) (\mathcal{S}^Q(E_{\delta})+\eps).
\end{align*}
then $\eps \rightarrow 0$ and $1-\delta<1$ implies
$\mathcal{S}^Q(E_{\delta})=0.$

\noindent
$(ii).$ Next we prove that
$$\mathcal{S}^{Q}(E_t)=0, \quad \forall\, t>1,$$ where
$$E_t:=\{q\in A: \mathcal{S}^{Q}(A\cap B(q,r))>t(2r)^Q, \forall\,  r \, \text{small enough}\}.$$

Now let $U$ be an open set such that $E_t \subset U$ and
$\mathcal{S}^Q(A\cap U)<\mathcal{S}^Q(E_t)+\eps$.
We define
$$\mathcal{F}:=\left\{ B(q,r): q\in E_t, B(q,r)\subset U, \tx{diam}\, B(q,r)\leq \delta \right\}.$$
Now we can apply Vitali covering lemma to $\mc{F}$ and get a family
$\{B_i\}$ of disjoint balls such that
$\mathcal{S}^{Q}(E_t \setminus \bigcup_i B_i)=0$.
Then we get
\begin{align*}
\mathcal{S}^{Q}(E_t)+\eps & > \mathcal{S}^{Q}(A \cap U)\\
                          &\geq  \mathcal{S}^{Q}(A\cap \bigcup B_i) \\
                          &\geq  t \sum \tx{diam}(B_i)^Q \\
                          &\geq t \, \mathcal{S}^Q_{\delta}(E_t \cap \bigcup B_i)\\
                          &\geq t \, \mathcal{S}^Q_{\delta}(E_t).
\end{align*}
Letting $\eps,\delta \rightarrow 0$ we have an absurd because
$t>1$.
\end{proof}
Since $A$ is open, from this lemma follows formula \eqref{eq}.

\brem
Notice that, for a $n$-dimensional Riemannian manifold, the tangent spaces at different points are isometric. As a consequence the Riemannian volume of the unit ball in the tangent space is constant and one can show that it is $C_n = \pi^\frac{n}{2}/ \Gamma(\frac{n}{2} + 1)$. Formula \eqref{eq:formulabella}, where $\mu=\tx{Vol}$ is the Riemannian volume, implies the well-known relation between $\tx{Vol}$ and the (spherical) Hausdorff measure
$$d\tx{Vol}\,=\,\frac{C_{n}}{2^{n}}\, d\mc{S}^{n}\,=\,\frac{C_{n}}{2^{n}}\, d\mc{H}^{n}.$$
\erem
    
\subsection{Proof of Corollary \ref{c-cont}: continuity of the density}


%

%

In this section we prove Corollary \ref{c-cont}. More precisely we study the continuity of the map 
\bqn \label{eq:1}f_{\mu\cal{S}}:q \mapsto \widehat{\mu}_{q}(\widehat{B}_{q}).\eqn 
To this purpose, it is sufficient to study the regularity under the hypothesis that $\widehat{\mu}_{q}$ does not depend on the point. Indeed it is easily seen that the smooth measure $\mu$, which is  defined on the manifold, induces on the nilpotent approximations a smooth family of measures $\{ \widehat{\mu}_{q}\}_{q\in M}$. In the case $\mu=\mc{P}$, this is a consequence of equality \eqref{eq:popp}. In other words we can identify all tangent spaces in coordinates with $\R^{n}$ and fix a measure $\widehat{\mu}$ on it.

We are then reduced to study the regularity of the volume of the unit ball of a smooth family of nilpotent structures in $\R^{n}$, with respect to a fixed smooth measure. Notice that this family depends on an $n$-dimensional parameter.

To sum up, we have left to study the regularity of the map
\bqn \label{reduced} q \mapsto \mc{L}(\widehat{B}_{q}),\qquad q\in M,
\eqn
where $\widehat{B}_{q}$ is the unit ball of a family of nilpotent structures $\widehat{\sub}_{q}$ in $\R^{n}$ and $\mc{L}$ is the standard Lebesgue measure.

Let us denote $\widehat{d}_{q}$ the sub-Riemannian distance in $\widehat{\sub}_{q}$ and $\rho_{q}:=\widehat{d}_{q}(0,\cdot)$. Following this notation $\widehat{B}_{q}=\{x\in \R^{n} | \, \rho_{q}(x)\leq 1\}$ and 
the coordinate expression \eqref{dil} implies that
\bqn \label{eq:111}
\mc{L}(\delta_{\alpha}(\widehat{B}_{q}))=\alpha^{Q} \mc{L}(\widehat{B}_{q}),\qquad \quad \all \alpha >0.
\eqn
Notice that, since our sub-Riemannian structure is regular, we can choose privileged coordinates $\psi_{q}:O_{q}\to \R^{n}$ smoothly with respect to $q$. 
Let now $q^{\prime}\neq q$, there exists $\alpha=\alpha(q,q^{\prime})$ such that (see Remark \ref{palla})
\bqn \label{eq:222}\delta_\frac{1}{\alpha}\widehat{B}_{q^{\prime}}\subset \widehat{B}_{q}\subset \delta_{\alpha}\widehat{B}_{q^{\prime}}.\eqn
Using \eqref{eq:111}, \eqref{eq:222} and monotonicity of the volume we get
$$\left(\frac{1}{\alpha^{Q}}-1\right)\mc{L}(\widehat{B}_{q^{\prime}})\leq \mc{L}(\widehat{B}_{q})-\mc{L}(\widehat{B}_{q^{\prime}})\leq (\alpha^{Q}-1)\mc{L}(\widehat{B}_{q^{\prime}}).$$
Then it is sufficient to show that $\alpha(q,q^{\prime})\to 1$ when $q^{\prime}\to q$. This property follows from the next 
\bl The family of functions $\rho_{q}|_{K}$ is equicontinuous for every compact $K\subset \R^{n}$. Moreover $\rho_{q^{\prime}}\to \rho_{q}$ uniformly on compacts in $\R^{n}$, as $q^{\prime}\to q$.
\el
In the case in which $\{\rho_{t}\}_{t>0}$ is the approximating family of  the nilpotent distance $\widehat{\rho}$,
this result is proved in \cite{agrachevjp}. See also \cite{nostrolibro} for a more detailed proof, using cronological calculus.
With the same arguments one can extend this result to any smooth family of regular sub-Riemannian structures. The key point is that we can construct a basis for the tangent space to the structure with bracket polynomials of the orthonormal frame where the structure of the brackets does not depend on the parameter.

\section{Proof of Theorem \ref{t-c4c5}: differentiability of the density in the corank 1 case}\label{s:diff}
In this section we prove Theorem \ref{t-c4c5},
We start by studying the contact case. Then we complete our analysis by reducing the quasi-contact case and the general case to the contact one.
\subsection{Normal form of the nilpotent contact case}\label{s-norm}

Consider a 2-step nilpotent sub-Riemannian manifold in $\R^{n}$ of rank $k$. 
 
 Select a basis $\{X_1,\ldots,X_k,Z_1,\ldots,Z_{n-k}\}$ such that 
\bqn\label{eq:sotto}
\begin{cases}
\displaystyle{\distr=\tx{span}\{X_{1},\ldots,X_{k}\} },\\
\displaystyle{[X_i,X_j]=\sum_{h=1}^{n-k}b_{ij}^h Z_h,  \qquad i,j=1,\ldots,k, \quad \tx{where} \quad b_{ij}^{h}=-b_{ji}^{h}},\\
[X_i,Z_j]=[Z_j,Z_h]=0, \qquad i=1,\ldots,k, \qquad j,h=1,\ldots,n-k. 
\end{cases}
\eqn
Hence the Lie bracket can be considered as a map
\bqn\label{eq:bra}[\cdot,\cdot]:\distr \times \distr \longrightarrow TM/\distr
\eqn
and is represented by the $n-k$ skew-simmetric matrices
$L^{h}=(b_{ij}^{h}), h=1,\ldots,n-k$.

In the contact case we have $(k,n)=(2\ell,2\ell+1)$ and our structure is represented by one non degenerate skew-symmetric matrix $L$. 
Take coordinates in such a way that $L$ is normalized in the following block-diagonal form

$$
L=\begin{pmatrix}
B_{1} & & \\
 & \ddots & \\
 && B_{\ell}
\end{pmatrix}
, \qquad \tx{where}\qquad
B_{i}:=\begin{pmatrix}
0 & -b_{i} \\
b_{i} & 0 \\
\end{pmatrix}, \qquad b_{i}>0.
$$
with eigenvalues $\pm i b_{1},\ldots,\pm ib_{\ell}$. Hence we can find vector fields $\{X_1,\ldots,X_\ell,Y_1,\ldots,Y_{\ell},Z\}$ such that relations \eqref{eq:sotto} reads

\bqn\label{eq:lie}
\begin{cases}
\distr=\tx{span}\{X_{1},\ldots,X_{\ell},Y_{1},\ldots,Y_{\ell}\},\\
[X_i,Y_i]=-b_i Z, \qquad~~~~~~~~~~~~~~~~~~~~~~ i=1,\ldots,\ell \\
[X_i,Y_j]=0,   \qquad~~~~~~~~ ~~~~~~~~~~~~~~~~~~~      i\neq j\\
[X_i,Z]=[Y_i,Z]=0,\qquad~~~~~~~~~~~~~~~~ i=1,\ldots, \ell 
\end{cases}
\eqn
In the following we call $b_1,\ldots,b_{\ell}$  \emph{frequences} of the contact structure.

We can recover the product on the group by the Campbell-Hausdorff formula. If we denote points $q=(x,y,z)$,
where
\begin{align*}
 x&=(x_1,\ldots,x_{\ell})\in \R^{\ell}, \qquad y=(y_1,\ldots,y_{\ell})\in \R^{\ell}, \qquad z\in \R,
\end{align*}
we can write the group law in coordinates 
\bqn \label{eq:coordpr}
q\cdot q^{\prime} = \left( x+x^{\prime},y+y^{\prime},z+z^{\prime}-\frac{1}{2}\sum_{i=1}^{\ell} b_i (x_ix^{\prime}_i-y_iy^{\prime}_i)\right).
\eqn
Finally, from \eqref{eq:coordpr}, we get the coordinate expression of the left-invariant vector fields of the Lie algebra, namely 
\begin{align} \label{eq:coordvf}
X_i&=\partial_{x_i}+\frac{1}{2} b_{i} y_i \partial_{z}, \qquad i=1,\ldots,\ell, \notag \\
Y_i&=\partial_{y_i}-\frac{1}{2} b_{i} x_i \partial_{z}, \qquad i=1,\ldots,\ell, \\
Z&=\partial_{z}. \notag
\end{align}
In this expression one of frequences $b_{i}$ can be normalized to 1.
\subsection{Exponential map in the nilpotent contact case}\label{s-exp}

Now we apply the PMP to find the exponential map $\Exp_{q_{0}}$ where $q_{0}$ is the origin. Define the 
hamiltonians (linear on fibers)
$$
h_{u_{i}}(\lam)=\la \lam, X_{i}(q)\ra, \quad~~~ 
h_{v_{i}}(\lam)=\la \lam, Y_{i}(q)\ra, \quad~~~
h_{w}(\lam)=\la \lam, Z(q)\ra. 
$$
Recall from Section \ref{s:pmp} that $q(t)$ is a normal extremal if and only if there exists $\lam(t)$ such that
\bqn\label{eq:hamiltonian}
\begin{cases}
\dot{u}_i=-b_i w v_{i}  \\
\dot{v}_{i}=b_iw u_{i}\\
\dot{w}=0 
\end{cases}
\qquad
\begin{cases}
\dot{x}_i=u_i \\
\dot{y}_i=v_i \\
\dot{z}=\frac{1}{2}\sum_{i}b_i(u_iy_{i}-v_{i}x_i)  
\end{cases}
\eqn
where
$$
u_{i}(t):=h_{u_{i}}(\lam(t)), \quad ~~~
v_{i}(t):=h_{v_{i}}(\lam(t)), \quad ~~~
w(t):=h_{w}(\lam(t)).
$$
\brem \label{rem:euclsr}
Notice that from \eqref{eq:hamiltonian} it follows that the sub-Riemannian length of a geodesic coincide with the Euclidean length of its projection on the horizontal subspace $(x_{1},\ldots,x_{n},y_{1},\ldots,y_{n})$.
$$l(\g)=\int_{0}^{T} \left( \sum_{i} (u_{i}^{2}(t)+v_{i}^{2}(t))\right)^{\frac{1}{2}} dt.$$
\erem
Now we solve \eqref{eq:hamiltonian} with initial conditions (see also Remark \ref{r-1/2})
\begin{gather}
(x^{0},y^{0},z^{0})=(0,0,0),\\
(u^0,v^0,w^0)=(u^0_1,\ldots,u^0_{\ell},v^0_1,\ldots,v^0_{\ell},w^0)\in S^{2\ell-1} \times \R .
\end{gather}
Notice that $w\equiv w^0$ is constant on geodesics. We consider separately the two cases:
\bi
\iii[$(i)$] If $w \neq 0$, we have (denoting $a_{i}:=b_{i}w$)
\begin{align}
u_i(t)&=u^0_i \cos a_i t - v^0_{i} \sin a_it, \notag \\
v_{i}(t)&=u^0_i \sin a_it + v^0_{i} \cos a_it,  \\
w(t)&=w. \notag
\end{align}
From \eqref{eq:hamiltonian}
one easily get
\begin{align} \label{eq:exp1}
x_i(t)&=\frac{1}{a_i}(u^0_i \sin a_it + v^0_{i} \cos a_it -v^0_{i}), \notag \\
y_{i}(t)&=\frac{1}{a_i}(-u^0_i \cos a_it + v^0_{i} \sin a_it +u^0_{i}), \\
z(t)&=\frac{1}{2w^2}(wt-\sum_i \frac{1}{b_i}  ((u^0_i)^2+(v^0_{i})^2))\sin a_it ). \notag
\end{align}

\iii[$(ii)$] If $w=0$, we find equations of straight lines on the horizontal plane in direction of the  vector $(u^0,v^0)$:
\begin{align*}
x_i(t)=u^0_it \qquad \quad
y_i(t)=v^0_it\qquad \quad
z(t)=0. 
\end{align*}
\ei
\brem To recover symmetry properties of the exponential map it is useful to rewrite \eqref{eq:exp1} in polar coordinates, using the following change of variables
\bqn \label{eq:change}
u_i= r_i \cos \theta_i, \qquad
v_i= r_i \sin \theta_i, \qquad 
 i=1,\ldots, \ell.
\eqn
In these new coordinates \eqref{eq:exp1} becomes
\begin{align} \label{eq:exp2}
x_i(t)&=\frac{r_i}{a_i}(\cos (a_it +\theta_i) - \cos \theta_i),\notag \\
y_{i}(t)&=\frac{r_i}{a_i}(\sin (a_it +\theta_i) - \sin\theta_i), \\
z(t)&=\frac{1}{2w^2}(wt-\sum_i \frac{r_i^2}{b_i} \sin a_it ). \notag
\end{align}
and the condition $(u^0,v^0)\in S^{2\ell-1}$ implies that $r=(r_1,\ldots,r_\ell)\in S^\ell$.

From equations \eqref{eq:exp2}  we easily see that the projection of a geodesic on every 2-plane $(x_i,y_i)$ is a circle, with period $T_i$, radius $\rho_i$ and center $C_i$ where 
\bqn \label{eq:circ}
T_{i}=\frac{2\pi}{b_i w}, \qquad\rho_i=\frac{r_i}{b_iw} \qquad C_i=-\frac{r_i}{b_iw}(\cos \theta_i,\sin \theta_i), \qquad \all i=1, \ldots,\ell
\eqn 
Moreover \eqref{eq:hamiltonian} shows that the $z$ component of the geodesic at time $t$ is the weighted sum (with coefficients $b_i$) of the areas spanned by the vectors  $(x_{i}(t),y_{i}(t))$ in  $\R^{2}$.
\erem

\bl \label{l:cutconj} Let $\g(t)$ be a geodesic starting from the origin and corresponding to the parameters $(r_i,\theta_i,w)$. The cut time $t_{c}$ along $\g$ is equal to first conjugate time and satisfies
\bqn \label{eq:conj}
t_c=\frac{2 \pi}{w \max_i b_i},
\eqn
with the understanding that $t_{c}=+\infty$, if $w=0$.
\el
\begin{proof} We divide the proof into two steps. Recall that a geodesic lose optimality either at a cut time or at a conjugate time, see Definition \ref{def:cutconj}. First we prove that \eqref{eq:conj} is a conjugate time and then that for every $t<t_c$ our geodesic is optimal. 

The case $w=0$ is trivial. Indeed the geodesic is a straight line and, by Remark \ref{rem:euclsr}, we have neither cut nor conjugate points. Then it is not restrictive to assume $w\neq 0$. Moreover, up to relabeling indices, we can also assume that $b_{1}\geq b_{2}\geq \ldots\geq b_{\ell}\geq 0$.\\
$(i)$ Since, by assumption, $b_{1}=\max_{i} b_i$, from \eqref{eq:exp2} it is easily seen that projection on the $(x_{1},y_{1})$-plane satisfies
$$x_{1}(t_c)=y_{1}(t_c)=0.$$
Next consider the one parametric family of geodesic with initial condition
$$(r_{1},r_{2},\ldots,r_{\ell},\theta_{1}+\phi,\theta_{2},\ldots,\theta_{\ell},w),\qquad \phi \in [0,2\pi].$$
It is easily seen from equation \eqref{eq:exp2} that all these curves have the same endpoints. Indeed neither $(x_{i},y_{i})$, for $i> 1$, nor $z$ depends on this variable. Then it follows that $t_{c}$ is a critical time for exponential map, hence a conjugate time.\\
$(ii)$ Since $w\neq 0$, our geodesic is non horizontal (i.e. $z(t)\not\equiv 0$). By symmetry, we can focus on the case $w>0$. We know that, for every $i$, the projection of every non horizontal geodesic on on the plane $(x_i,y_i)$ is a circle. Moreover for the $i$-th projected curve, the distance from the origin is easily computed
$$\eta_{i}(t)=\sqrt{x_{i}(t)^{2}+y_{i}(t)^{2}}=r_{i}\, t \, \sin_{c}(\frac{b_i w t}{2}), \qquad \text{where}\qquad \sin_{c}(x)=\frac{\sin x}{x}.$$
Let now $\bar t<t_c$, we want to show that there are no others geodesics $\til{\g}(t)$, starting from the origin,  that reach optimally the point 
$\g(\bar t)$ at the same time $ \bar t$. 
Assume that $\til{\g}(t)$ is associated to the parameters $(\til{r}_{i},\til{\theta}_{i},\til{w})$, where $(\til{r_{1}},\ldots,\til{r}_{\ell})\in S^{\ell}$, and let us argue by contradiction. 
If $\g(\bar t)=\til{\g}(\bar t)$ it follows that $\eta_{i}(\bar t)=\til{\eta}_{i}(\bar t)$ for every $i$, that means
\bqn \label{eq:sds}
r_{i}\, \bar{t}\sin_{c}(\frac{b_i w\bar t}{2})=\til{r}_{i}\, \bar{t} \, \sin_{c}(\frac{b_i \til{w} \bar t}{2}),\qquad i=1,\ldots,\ell.
\eqn
Notice that, once $\til{w}$ is fixed, $\til{r}_{i}$ are uniquely determined by \eqref{eq:sds} (recall that $\bar t$ is fixed). Moreover, $\til{\theta}_{i}$ also are uniquely determined by relations \eqref{eq:circ}. Finally, from the assumption that $\til{\g}$ also reach optimally the point $\til{\g}(\bar t)$, it follows that
\bqn \label{aaaaa}
\bar t<t_{c}=\frac{2\pi}{b_{1}\til{w}} \qquad \Longrightarrow \qquad \frac{b_i \til{w} \bar t}{2}< \pi \quad \all i=1,\ldots,\ell.
\eqn
Since $\sin_{c}(x)$ is a strictly decreasing function on $[0,\pi]$ it follows from \eqref{eq:sds} that, if $\til{w}\neq w$ we have
\begin{gather*}
\til{w}>w\ \  \Rightarrow \ \ \til{r}_{i}>r_{i} \ \  \all i=1,\ldots,\ell \ \ \Rightarrow \ \ \sum_{i} \til{r}_{i}^{2}>\sum_{i} r_{i}^{2}=1,\\
\til{w}<w\ \  \Rightarrow \ \ \til{r}_{i}<r_{i} \ \  \all i=1,\ldots,\ell \ \ \Rightarrow \ \ \sum_{i} \til{r}_{i}^{2}<\sum_{i} r_{i}^{2}=1,
\end{gather*}
which contradicts the fact that $(\til{r_{1}},\ldots,\til{r}_{\ell})\in S^{\ell}$.

\end{proof}

Consider now the exponential map from the origin
 \bqn \label{eq:exp3}
(r,\theta, w)\mapsto \mathsf{Exp}(r,\theta,w)=
\begin{cases} 
x_i=\dfrac{r_i}{b_iw}(\cos (b_iw +\theta_i) - \cos \theta_i) \\[0.4cm]
y_{i}=\dfrac{r_i}{b_iw}(\sin (b_iw +\theta_i) - \sin \theta_i) \\[0.2cm]
z=\dfrac{1}{2w^2}(w|r|^2-\sum_i \dfrac{r_i^2}{b_i} \sin b_iw )=\sum_i\dfrac{r_i^2}{2b_iw^2}(b_iw- \sin b_iw ).  
\end{cases}
\eqn

By Lemma \ref{l:cutconj}, the set $D$ where geodesics are optimal and their length is less or equal to 1, is characterized as follows
$$
D:=\{(r,\theta,w), |r|\leq 1, |w|\leq 2\pi/\max b_i\}.$$
Thus the restriction of the exponential map to the interior of $D$ gives a regular parametrization of the nilpotent unit ball $\widehat{B}$, and we can compute its volume with the change of variables formula
\bqn
\mc{L}(\widehat{B})=\int_{\widehat{B}}d\mc{L}=\int_{D} |\det J |\, R \, drd\theta dw,
\eqn
where $J$ denotes the Jacobian matrix of the exponential map \eqref{eq:exp3}. Notice that we have to integrate with respect to the measure 
$dudvdw=R \, drd\theta dw$, where $R=\prod_i r_i$
because of the change of variables \eqref{eq:change}.

\bl The Jacobian of Exponential map \eqref{eq:exp3} is given by the formula
\bqn \label{eq:jac2}
\det J(r,\theta,w)=
\frac{4^\ell R}{B^2 w^{2\ell+2}}\sum_{i=1}^\ell \left( \prod_{j\neq i}\sin^2(\frac{b_jw}{2})\right) \sin(\frac{b_iw}{2})\left(	\frac{b_iw}{2}\cos(\frac{b_iw}{2})-\sin(\frac{b_iw}{2})\right)r_i^2,\eqn
where 
we denote with 
$B=\prod_i b_i$.
\el
\begin{proof}
We reorder variables in the following way 
$$(r_1,\theta_1,\ldots,r_\ell,\theta_\ell,w), \qquad (x_1,y_1,\ldots,x_\ell,y_\ell,z),$$
in such a way that the Jacobian matrix $J$ of the exponential map \eqref{eq:exp3} is
\bqn \label{eq:J}
J=
\begin{pmatrix}
Q_1&      &            &      & W_1 \\
      &Q_2&            &      & W_2\\
      &       &\ddots&      & \vdots\\
      &       &           &Q_\ell &W_\ell\\
Z_1 & Z_2& \ldots & Z_\ell & \partial_wz    
\end{pmatrix}
\eqn
where we denote
\begin{align}
Q_i&=\begin{pmatrix} Q_i^r & Q_i^{\theta}\end{pmatrix}:=\begin{pmatrix}
\de{x_i}{r_i}& \de{x_i}{\theta_i}  \\
\de{y_i}{r_i}& \de{y_i}{\theta_i}
\end{pmatrix},  \qquad i=1,\ldots,\ell, \\
 W_i&:=\begin{pmatrix} \de{x_i}{w}\\ \de{y_i}{w}\end{pmatrix}, \qquad Z_i:=\begin{pmatrix} \de{z}{r_i}&\de{z}{\theta_i}\end{pmatrix}. \notag
 \end{align}
Notice that $x_i,y_i$ depend only on $r_i,\theta_i$.

To compute the determinant of $J$, we write $$z=z(r,w)=\sum_i z_i(r_i,w), \qquad z_i:=\dfrac{r_i^2}{2b_iw^2}(b_iw- \sin b_iw ),$$
and we split the last column of $J$ as a sum
\bqn \label{eq:x1}\begin{pmatrix}
W_1\\
W_2\\
\vdots\\
W_\ell \\
\partial_w z
\end{pmatrix}=
\begin{pmatrix}
W_1\\
0\\
\vdots\\
0\\
\partial_w z_1
\end{pmatrix}+
\begin{pmatrix}
0\\
W_2\\
\vdots\\
0\\
\partial_w z_2
\end{pmatrix}+\cdots+
\begin{pmatrix}
0\\
0\\
\vdots\\
W_\ell\\
\partial_w z_\ell
\end{pmatrix}.
\eqn 
Notice that in the $i$-th column only the $i$-th variables appear.

By multilinearity of determinant, $\det J$ is the sum of the determinants of $\ell$ matrices, obtained by replacing each time last column with one of vectors appearing in the sum \eqref{eq:x1}. If we replace it, for instance, with the first term, we get
$$
J_1=\begin{pmatrix}
Q_1&      &            &      & W_1 \\
      &Q_2&            &      & 0\\
      &       &\ddots&      & \vdots\\
      &       &           &Q_\ell&0\\
Z_1 & Z_2& \ldots & Z_\ell & \partial_w z_1\\      
\end{pmatrix}
$$
Now with straightforward computations (notice that $\partial_{\theta}z_{1}=0$ in $Z_{1}$), we get
\begin{align*}
\det J_1&= \det Q_2 \cdots \det Q_\ell \cdot
\det \begin{pmatrix}
Q_1 & W_1\\
Z_{1}& \de{z_1}{w}
\end{pmatrix}
\\
&=\det Q_2 \cdots \det Q_\ell \cdot
\det \begin{pmatrix}
Q_1^r & Q_1^{\theta} & W_1\\
\de{z_1}{r_1} & 0 & \de{z_1}{w}
\end{pmatrix}.
\end{align*}
Setting
$$A_i=
\begin{pmatrix} 
Q_i^{\theta} &W_i\end{pmatrix}\qquad i=1,\ldots,\ell,$$
we find
$$\det J_1=\det Q_2 \cdots \det Q_\ell \cdot (\partial_w z_1 \det Q_1 + \partial_{r_1} z_1 \det A_1).$$
Similarly, we find analogous expressions for $J_2,\ldots,J_\ell$. Then
$$\det J= \sum_{i=1}^{\ell} \det J_i= \sum_{i=1}^{\ell} (\prod_{j\neq i} \det Q_j)(\partial_w z_i \det Q_i + \partial_{r_i} z_i \det A_i).$$
From \eqref{eq:exp3}, by direct computations, it follows
\begin{align*}
\det Q_j&= \frac{4r_j}{b_j^2w^2} \sin^2(\frac{b_jw}{2}) ,
\end{align*}
Moreover, with some computations, one can get
\begin{align*}
\partial_w z_i&=\frac{r_{i}^{2}}{w^{2}}\sin^{2}(\frac{b_{i}w}{2})-\frac{r_{i}^{2}}{b_{i}w^{3}}(b_iw-\sin b_{i}w)\\
\partial_{r_i} z_i&= \frac{r_i}{b_iw^2}(b_iw-\sin b_{i}w)\\
\det A_i&=-\frac{4r_i^2}{b_i^2w^3}\sin(\frac{b_iw}{2})\left(\frac{b_iw}{2} \cos(\frac{b_iw}{2})-\sin (\frac{b_iw}{2})\right)
\end{align*}
and finally, after some simplifications
\begin{align*}
\partial_w z_i \det Q_i + \partial_{r_i} z_i \det A_i&= -\frac{4r_i^3}{b_i^2w^4}\sin(\frac{b_iw}{2})\left(\frac{b_iw}{2} \cos(\frac{b_iw}{2})-\sin (\frac{b_iw}{2})\right).
\end{align*}
from which we get \eqref{eq:jac2}. 
\end{proof}

From the explicit expression of the Jacobian \eqref{eq:jac2} we see that integration with respect to horizontal variables $(r_{i},\theta_{i})$ does not involve frequences, providing a costant $C_{\ell}$. Hence, the computation of the volume reduces to a one dimensional integral in the vertical variable $w$:
$$V=\int_{-\frac{2\pi}{\max b_i}}^{\frac{2\pi}{\max b_i}}\frac{C_\ell}{B^{2}w^{2\ell+2}}\sum_{i=1}^\ell \left( \prod_{j\neq i}\sin^2(\frac{b_jw}{2})\right) \sin(\frac{b_iw}{2})\left(\frac{b_iw}{2}\cos(\frac{b_iw}{2})-\sin(\frac{b_iw}{2})\right)dw.$$
Using simmetry property of the jacobian with respect to $w$ and making the change of variable $2s=w$, the volume become (reabsorbing all costants in $C_{\ell}$)
\begin{equation} \label{eq:drop}
V=\int_0^{\frac{\pi}{\max b_i}} \frac{C_\ell}{ B^{2}s^{2\ell+2}}\sum_{i=1}^\ell \left( \prod_{j\neq i}\sin^2(b_js)\right) \sin(b_is)(b_is\cos(b_is)-\sin(b_is)) ds.
\end{equation}

\brem In the case $\ell=1$ (Heisenberg group) we get
$$V=\frac{1}{12}(1+2\pi \,\text{Si}(2\pi)) \simeq 0.8258, \qquad \text{Si}(x):=\int_{0}^{x}\frac{\sin t}{t}dt.$$
Notice that this is precisely the value of the constant $f_{\mc{PS}}$, since in this case Popp's measure coincide with the Lebesgue measure in our coordinates.
\erem

\subsection{Differentiability properties: contact case}


Let us come back to the differentiability of the map 
\bqn\label{eq:npar} 
q \mapsto \mc{L}(\widehat{B}_{q}).
\eqn

A smooth family of sub-Riemannian structures is represented by a smooth familly of skew symmetric matrices $L(q)$ (see Section \ref{s-norm}).
Recall that, for a smooth family of skew-symmetric matrices that depend on a $n$-dimensional parameter, the eigenvalue functions $q\mapsto b_{i}(q)$ exists and are Lipschitz continuous with respect to $q$ (see \cite{kato}). 

Thus, if we denote with $V(q)$ the volume of the nilpotent unit ball corresponding to frequences $b_{1}(q),\ldots, b_{\ell}(q)$, formula \eqref{eq:drop} can be rewritten as
\begin{equation} \label{eq:volume}
V(q)=\int_0^{a(q)} G(q,s) ds,
\end{equation}
where 
\begin{gather} 
a(q):=\frac{\pi}{\max b_{i}(q)}, \notag\\
G(q,s):=\frac{1}{ s^{2\ell+2}B^{2}}\sum_{i=1}^\ell \left( \prod_{j\neq i}\sin^2(b_j(q)s)\right) \sin(b_i(q)s)(b_i(q)s\cos(b_i(q)s)-\sin(b_i(q)s)).\label{eq:explicit}
\end{gather}

Notice that we have dropped the constant $C_{\ell}$ that appear in \eqref{eq:drop} since it does not affect differentiability of the volume.

\brem
Since the family of sub-Riemannian structures $q\mapsto L(q)$ is smooth, the exponential map smoothly depends on the point $q$. As a consequence the integrand $G(q,s)$, being the Jacobian of the exponential map, is a smooth function of its variables. 

In addition, altough $q\to b_{i}(q)$ is only Lipschitz, it is easy to see that the function $a(q)= \max b_{i}(q)$ is semiconvex with respect to $q$ (see also \cite{semicbook}). In particular $a(q)$ admits second derivative almost everywhere.
\erem

If for all $q$ all $b_{i}(q)$ are different (i.e. there are no resonance points), then the eigenvalue functions can be chosen in a smooth way. As a consequence the volume $V(q)$ is smooth, since all functions that appear in \eqref{eq:volume} are smooth. This argument provides a proof of Proposition \ref{p} in the case of contact structures.

On the other hand, we prove that, along a curve where $a(q)$ is not smooth, i.e. when the two bigger frequences cross, $V(q)$ is no longer smooth at that point.
\begin{proof}[Proof of Theorem \ref{t-c4c5}] 

In this proof by a resonance point we mean a point $q_{0}$ 
where  the two (or more) biggest frequences of $L(q_{0})$ coincide. More precisely, if we order the frequences as
$$b_{1}(q)\geq b_{2}(q) \geq \ldots \geq b_{\ell}(q),$$
a resonance point is defined by the relation $b_{1}(q_{0})=b_{2}(q_{0})$.
As we noticed, $V$ is smooth at non-resonance points.

We divide the proof into two steps: first we prove that $V\in\con^{3}$ and then we show that, in general, it is not smooth (but $\con^{4}$ on smooth curves).

$(i)$. We have to show that $V$ is $\con^{3}$ in a neighborhood of every resonance point $q_{0}$. First split the volume as follows
\bqn \label{eq:VW}
V(q)=\int_{0}^{a(q_{0})}G(q,s) ds+\int_{a(q_{0})}^{a(q)}G(q,s) ds
\eqn
The first term in the sum is smooth with respect to $q$, since it is the integral of a smooth function on a domain of integration that does not depend on $q$. We are then reduced to the regularity of the function
\bqn \label{eq:W}
W(q):=\int_{a(q_{0})}^{a(q)}G(q,s) ds
\eqn

We have the following key estimate
\bl \label{l:asda}
Let $q_{0}\in M$ be a resonance point. Then, for any neighborhood of $q_{0}$, there exists $C>0$ such that
\bqn \label{eq:stima}
\left|\int_{a(q_{0})}^{a(q)} G(q,s) ds\right|\leq C|q-q_{0}|^{4}
\eqn
\el
\begin{proof} It is sufficient to prove that every derivative up to second order of $G$ vanish at $(q_{0},a(q_{0}))$. Indeed, being $G$ smooth, computing its Taylor polynomial at $(q_{0},a(q_{0}))$ only terms with order greater or equal than three appear (both in $q-q_{0}$ and $s-a(q_{0})$). Thus, integrating with respect to $s$ and using that $|a(q)-a(q_{0})|= O(|q-q_{0}|)$, we have the desired result.

From the explicit formula \eqref{eq:explicit} it is easy to see that $G(q,a(q))\equiv0$ for every $q\in M$. In particular $G(q_{0},a(q_{0}))=0$.
Moreover, 
since at a resonance point $q_{0}$ at least the two bigger eigenvalues coincide, say $b_{1}(q_{0})=b_{2}(q_{0})=\beta$, we have $a(q_{0})=\pi/\beta$ and in a neighborhood of $(q,s)=(q_{0},a(q_{0}))$
\bqn \label{eq:sin3}
\sin^{2}(b_{1}(q)\frac{\pi}{a(q)})\sin(b_{2}(q)\frac{\pi}{a(q)})=O(|b_{1}(q)-b_{2}(q)|^{3})=O(|q-q_{0}|^{3}).
\eqn
due to the Lipschitz property of $b_{j}(q)$, for $j=1,2$. 

From \eqref{eq:explicit} and \eqref{eq:sin3} one can easily get that every derivative of $G$ up to second order (in both variables $q$ and $s$) vanish at $(q_{0},a(q_{0}))$.
%

\end{proof}
To show that $V\in \con^{3}$ we compute the first three derivatives of $W$ at a non resonant points $q$, and we show that, when $q$ tend to a resonance point $q_{0}$, they tends to zero. We then conclude the continuity of the first three derivatives by Lemma \ref{l:asda}.
  
In the following, for simplicity of the notation, we will denote by $\frac{\partial}{ \partial q}$ the partial derivative with respect to some coordinate function on $M$. For instance $\frac{\partial^{2} W}{ \partial q^{2}}$ denote some second order derivative $\frac{\partial^{2} W}{ \partial x_{i} \partial x_{j}}$.

At non resonance points $q$ we have
\begin{align*}\frac{\partial W}{\partial q}(q)&= \underbrace{G(q,a(q))}_{=0}\frac{\partial a}{\partial q}(q)+ \int_{a(q_{0})}^{a(q)} \frac{\partial G}{\partial q}(q,s)ds \\
&=\int_{a(q_{0})}^{a(q)} \frac{\partial G}{\partial q}(q,s)ds\end{align*}
which tends to zero for $q \to q_{0}$. Using Lemma \ref{l:asda} we conclude $V\in \con^{1}$. Next let us compute
 
\bqn \label{eq:sss}
\frac{\partial^{2} W}{\partial q^{2}}(q)=\frac{\partial G}{\partial q}(q,a(q))\frac{\partial a}{\partial q}(q)+ \int_{a(q_{0})}^{a(q)} \frac{\partial^{2} G}{\partial q^{2}}(q,s)ds
\eqn
Since $\frac{\partial G}{\partial q}(q_{0},a(q_{0}))=0$ (see proof of Lemma \ref{l:asda}) and $\frac{\partial a}{\partial q}$ is bounded by Lipschitz continuity of the maximum eigenvalue, it follows that $\frac{\partial^{2} W}{\partial q^{2}}$ tends to zero as $q\to q_{0}$, and using again Lemma \ref{l:asda}, we have that $V\in \con^{2}$. 
In analogous way one can compute the third derivative
\begin{align} \label{eq:sss}
\frac{\partial^{3} W}{\partial q^{3}}(q)=\frac{\partial^{2} G}{\partial q \partial s}&(q,a(q))\left(\frac{\partial a}{\partial q}(q)\right)^{2}+2\frac{\partial^{2} G}{\partial q^{2}}(q,a(q))\frac{\partial a}{\partial q}(q)+\nn\\
&+\frac{\partial G}{\partial q}(q,a(q))\frac{\partial^{2} a}{\partial q^{2}}(q)+ \int_{a(q_{0})}^{a(q)} \frac{\partial^{3} G}{\partial q^{3}}(q,s)ds
\end{align}
Using again that every second derivative of $G$ vanish at $(q_{0},a(q_{0}))$ and that $\frac{\partial a}{\partial q}$ is bounded, it remains to check that $\frac{\partial G}{\partial q}(q,a(q))\frac{\partial^{2} a}{\partial q^{2}}(q)$ tends to zero as $q \to q_{0}$. From \eqref{eq:sin3} one can see that $\frac{\partial G}{\partial q}=O(|b_{1}(q)-b_{2}(q)|^{2})$. Hence it is sufficient to prove that $\frac{\partial^{2} a}{\partial q^{2}}=O(1/|b_{1}(q)-b_{2}(q)|)$, which is a consequence of the following lemma.

\smallskip
Notice that, for every skew-symmetric matrix $A$, with eigenvalues $\pm i \lam_{j}$, with  $j=1,\ldots,n$, the matrix $iA$ is an Hermitian matrix which has eigenvalues $\pm \lam_{j},\, j=1,\ldots,n$.

\bl Let $A,B$ be two $n\times n$ Hermitian matrices and assume that for every $t$ the matrix $A+tB$ has a simple eigenvalue $\lam_{j}(t)$.
Then the the following equation is satisfied
\bqn \label{eq:lemmone}
\ddot{\lam}_{j}=2\sum_{k\neq j}\frac{|\langle B x_{j},x_{k}\rangle|^{2}}{\lam_{j}-\lam_{k}}
\eqn
where $\{x_{k}(t)\}_{k=1,\ldots,n}$ is an orthonormal basis of eigenvectors and $x_{j}(t)$ is the eigenvector associated to $\lam_{j}(t)$.
\el
\begin{proof} In this proof we endow $\C^{n}$ with the standard scalar product $\la z ,w\ra=\sum_{k=1}^{n} z_{k} \overline{w_{k}}$.
Since $\lam_{j}(t)$ is simple for every $t$, both $\lam_{j}(t)$ and the associated eigenvector $x_{j}(t)$ can be choosen smoothly with respect to $t$. By definition
$$(A+tB)x_{j}(t)=\lam_{j}(t)x_{j}(t), \qquad |x_{j}(t)|=1.$$
Then we compute the derivative with respect to $t$ of both sides
\bqn\label{eq:ausil}(A+tB)\dot{x}_{j}(t)+Bx_{j}(t)=\dot{\lam}_{j}(t)x_{j}(t)+\lam_{j}(t)\dot{x}_{j}(t),\eqn
and computing the scalar product with $x_{j}(t)$ we get
\bqn \label{eq:www}
\dot{\lam}_{j}(t)=\la Bx_{j}(t),x_{j}(t)\ra, \qquad \tx{hence} \qquad \ddot{\lam}_{j}(t)=2\, \tx{Re} \la Bx_{j}(t), \dot{x}_{j}(t)\ra.
\eqn
using that $A+tB$ is Hermitian and $\la\dot{x}_{k}(t),x_{k}(t)\ra=0$.
On the other hand, the scalar product of \eqref{eq:ausil} with $x_{k}$, with $k \neq j$ gives
$$\la\dot{x}_{j},x_{k}\ra=\frac{\la Bx_{j},x_{k}\ra}{\lam_{j}-\lam_{k}}.$$
Substituting $\dot{x}_{j}=\sum_{k= 1}^{n} \la\dot{x}_{j},x_{k}\ra x_{k}$ in \eqref{eq:www} we have \eqref{eq:lemmone}.
\end{proof}
Then the third derivative goes to zero for $q\to q_{0}$, and using again Lemma \ref{l:asda} we conclude that $V\in \con^{3}$.

\medskip
$(ii)$. Now we study the restriction of the map \eqref{eq:npar} along any smooth curve on the manifold, and we see that, due to a simmetry property, $V$ is $\con^{4}$ on every curve but in general is not $\con^{5}$. Notice that \eqref{eq:stima} gives no information on the fourth derivative but the fact that it is bounded. Indeed it happens that it is continuous on every smooth curve but its value depend on the curve we choose. 

From now on, we are left to consider a smooth \emph{one-parametric} family of sub-Riemannian structure, i.e. of skew symmetric matrices.

\brem 
An analytic family of skew-simmetric matrices $t\mapsto L(t)$ depending on one parameter, can be simultaneously diagonalized (see again \cite{kato}), in the sense that there exists an analytic (with respect to $t$) family of orthogonal changes of coordinates and analytic functions $b_{i}(t)>0$ such that
\bqn \label{eq:Lt1}
L=\begin{pmatrix}
B_{1}(t) & & \\
 & \ddots & \\
 && B_{\ell}(t)
\end{pmatrix}
, \qquad \tx{where}\qquad
B_{i}(t):=\begin{pmatrix}
0 & -b_{i}(t) \\
b_{i}(t) & 0 \\
\end{pmatrix}.
\eqn
In the case of a $\con^{\infty}$ family $t\mapsto L(t)$, we can apply the previuos result to the Taylor polynomial of this family. As a consequence we get an approximate diagonalization  for $L(t)$, i.e. for every $N>0$ there exists a smooth family of orthogonal changes of coordinates and smooth functions $b_{i}(t)>0$ such that every entry out of the diagonal in $L(t)$ is $o(t^{N})$. Namely
\bqn \label{eq:Lt2}
L(t)=\begin{pmatrix}
B_{1}(t) & & o(t^{N}) \\
 & \ddots & \\
 o(t^{N})&& B_{\ell}(t)
\end{pmatrix}
, \qquad \tx{where}\qquad
B_{i}(t):=\begin{pmatrix}
o(t^{N}) & -b_{i}(t) \\
b_{i}(t) & o(t^{N}) \\
\end{pmatrix}.
\eqn

Since we are interested, in the study the $\con^k$ regularity  of \eqref{eq:npar}, for $k$ finite, in what follows we can ignore higher order terms and assume that $L(t)$ can be diagonalized as in the analytic case \eqref{eq:Lt1}.
\erem

\medskip
From the general analysis we know that  $V$ is $\con^{3}$. To prove that $t\mapsto V(t)$ is actually $\con^{4}$ 
we discuss first the easiest case $\ell=2$ (i.e. the contact (4,5) case) and then generalize to any $\ell$.\\

$(i)$ Case $\ell=2$. To start, assume that $b_{1}(t)$, $b_{2}(t)$ cross transversally at $t=0$ .
This means that for the volume $V(t)$ we have the expression
\bqn
V(t)=
\begin{cases}
\displaystyle{\int_{0}^{\frac{\pi}{b_{1}(t)}} G(t,s)\, ds}, &  \text{if }t>0,\\[0.3cm]
\displaystyle{\int_{0}^{\frac{\pi}{b_{2}(t)}} G(t,s)\, ds}, &  \text{if }t<0,\\
\end{cases} \quad  \text{ and } \quad 
\begin{matrix}
b_{1}(0)=b_{2}(0)\\[0.2cm]
b_{1}^{\prime}(0)\neq b_{2}^{\prime}(0)\\
\end{matrix}
\eqn


Since the regularity of the volume does not depend on the value $b_{1}(0)=b_{2}(0)$, we can make the additional assumption
$$b_{i}(t)=1+t \,c_{i}(t), \qquad i=1,2$$
for some suitable functions $c_{1}(t),c_{2}(t)$. Notice that $a^{\prime}(t)$ is discontinuous at $t=0$ and the left and right limits are
$$a^{\prime}_{+}:=\lim_{t\to 0+}a^{\prime}(t)=-\pi c_{1}(0), \qquad a^{\prime}_{-}:=\lim_{t\to 0-}a^{\prime}(t)=-\pi c_{2}(0),$$

 From the explicit expression of $G$ it is easy to compute that 
  \begin{gather*} \label{eq:ugconj2}
 \frac{\partial^{3} G}{\partial t^{3}}(0,a(0))=\frac{6}{\pi^{2}}c_{1}c_{2}(c_{1}+c_{2}), \\
  \frac{\partial^{3} G}{\partial t^{2} \partial s}(0,a(0))=\frac{2}{\pi^{3}}(c_{1}^{2}+4c_{1}c_{2}+c_{2}^{2}),  \qquad \frac{\partial^{3} G}{\partial t \partial s^{2}}(0,a(0))=\frac{6}{\pi^{4}}(c_{1}+c_{2}), 
\end{gather*}
where we denote for simplicity $c_{i}:=c_{i}(0)$.

Let us compute $\frac{\partial^{4} W}{\partial t^{4}}$. To this purpose let us differentiate with respect to $t$ formula \eqref{eq:sss} (where $q$ is replaced by $t$). Using the fact that all second derivatives of $G$ vanish at $(t,s)=(0,a(0))$ (see the proof of Lemma \ref{l:asda}) we have that the $4$-th derivative of $W$ at $t=0$ is computed as follows
\begin{align*}
\lim_{t\to 0+}W^{(4)}(t)&= 3 \frac{\partial^{3} G}{\partial t^{3}}a_{+}^{\prime}+ 3\frac{\partial^{3} G}{\partial t^{2} \partial s}(a_{+}^{\prime})^{2}+ \frac{\partial^{3} G}{\partial t \partial s^{2}}(a_{+}^{\prime})^{3}
,\\
\lim_{t\to 0-}W^{(4)}(t)&= 3 \frac{\partial^{3} G}{\partial t^{3}}a_{-}^{\prime}+ 3\frac{\partial^{3} G}{\partial t^{2} \partial s}(a_
{-}^{\prime})^{2}+ \frac{\partial^{3} G}{\partial t \partial s^{2}}(a_{-}^{\prime})^{3},
\end{align*}
where $W$ is defined in \eqref{eq:W}. It is easily checked that $W^{(4)}$ is continuous (but does not vanish!). Indeed we have
$$\lim_{t\to 0+}W^{(4)}(t)=\lim_{t\to 0-}W^{(4)}(t)=-\frac{12}{\pi}c_{1}^{2}c_{2}^{2}.
$$
The same argument produce an example that, in general, $V(t)$ is not $\con^{5}$. Assuming
$$b_{i}(t)=1+t \,c_{i}, \qquad c_{1}\neq c_{2} \quad \tx{constant}, \qquad i=1,2,$$
a longer computation, but similar to the one above, shows that
\begin{align*}
\lim_{t\to 0+}W^{(5)}(t)&= -\frac{2}{\pi} \, c_{1}^{3} (13c_{1}^{2}-29 c_{1}c_{2}+22 c_{2}^{2})
\\
\lim_{t\to 0-}W^{(5)}(t)&=-\frac{2}{\pi} \, c_{2}^{3} (13c_{2}^{2}-29 c_{1}c_{2}+22 c_{1}^{2})
\end{align*}
and the $5$-th derivatives do not coincide.
\brem
Notice that the assumption of transversality on $b_{1}(t)$ and $b_{2}(t)$ at $t=0$ is not restrictive. Indeed if $b_{1}(t)-b_{2}(t)=O(t^{k})$ for some $k>1$, then from the proof of Lemma \ref{l:asda} (see in particular \eqref{eq:sin3}) it follows that at least $2k+1$ derivatives of $G$ vanish at $(t,s)=(0,a(0))$, increasing the regularity of $W$. 
\erem

$(ii)$ General case. We reduce to case $(i)$. 

We can write $G(t,s)=\sum_{i=1}^{\ell} G_{i}(t,s)$ and $V(t)=\sum_{i=1}^{\ell} V_{i}(t)$ where we set
\begin{gather}
G_{i}(t,s):=\frac{1}{ s^{2\ell+2}} \left( \prod_{j\neq i}\sin^2(b_j(t)s)\right) \sin(b_i(t)s)(b_i(t)s\cos(b_i(t)s)-\sin(b_i(t)s)), \label{eq:gi}\\
V_{i}(t):=\int_0^{a(t)} G_{i}(t,s) ds, \qquad\qquad i=1,\ldots, \ell. \notag
\end{gather} 

Assume that $b_{1},b_{2}$ are the bigger frequences and that they cross at $t=0$, i.e.
\begin{gather*}
b_{i}(t)< b_{2}(t)< b_{1}(t), \qquad \all t<0, \quad \all i=3,\ldots,n.\\
b_{i}(t)<b_{1}(t)< b_{2}(t), \qquad \all t>0,  \quad \all i=3,\ldots,n.
\end{gather*}

From the explicit expression above it is easy to recognise that for $G_{1}$ and $G_{2}$  we can repeat the same argument used in $(i)$. Indeed if we denote with $\til{G}(t,s)$ the integrand of the $(4,5)$ case we can write $G_{1}+G_{2}$ as the product of a smooth function and $\til{G}$
$$G_{1}(t,s)+G_{2}(t,s)=\left( \frac{1}{ s^{2\ell-4}}  \prod_{j=3}^{\ell}\sin^2(b_j(t)s)\right)   \til{G}(t,s),$$
which implies that that $V_{1}+V_{2}$ is a $\con^{4}$ function.

Moreover it is also easy to see that $V_{3},\ldots,V_{n}$ are $\con^{4}$. Indeed from the fact that $b_{1}(t)$ and $b_{2}(t)$ both appear in $\sin^{2}$ terms in $G_{i}(t,s)$ for $i>3$, it follows that in this case
$$G_{i}(t,a(t))\equiv \frac{\partial G_{i}}{\partial t}(t,a(t))\equiv 0, \quad \frac{\partial^{2} G_{i}}{\partial t^{2}}(0,a(0))=0, \quad\frac{\partial^{3} G_{i}}{\partial t^{3}}(0,a(0))=0, \qquad  i=3,\ldots,n,$$ 
and we can apply the same argument used in $(i)$ to the function $V_{i}^{\prime}(t)=\int_{0}^{a(t)} \frac{\partial G_{i}}{\partial t}(t,s)ds$, for $i=3,\ldots,n$, showing that it is $\con^{3}$, that means $V_{i}\in \con^{4}$. 
 \end{proof}
\brem As we said the value of the 4-th derivative depend on the curve we chose, hence we cannot conclude that $V\in \con^{4}$ in general. Moreover, we explicitly proved that $V\notin \con^{\infty}$ since in general is not $\con^{5}$, even when restricted on smooth curves.

Moreover from the proof it also follows that, if more than two frequences coincide at some point (for instance if we get a triple eigenvalue), we have a higher order regularity for every $V_{i}$, and the regularity of $V$ increases.

\erem

\subsection{Extension to the quasi-contact case}
Recall that in the quasi contact case the dimension of the distribution is odd and the kernel of the contact form is one dimensional. Hence, applying the same argument used in Section \ref{s-norm},  we can always normalize the matrix $L$ in the following form:
$$L=\begin{pmatrix}
B_{1} & & &\\
 & \ddots & & \\
 && B_{\ell} &\\
 &&&0
\end{pmatrix}
, \qquad
B_{i}:=\begin{pmatrix}
0 & -b_{i} \\
b_{i} & 0 \\
\end{pmatrix}, \quad b_{i}>0.$$

In other words we can select a basis $\{X_{1},\ldots,X_{\ell},Y_{1},\ldots,Y_{\ell},K,Z\}$ such that
\bqn\label{eq:lie2}
\begin{cases}
\distr=\tx{span}\{X_{1},\ldots,X_{\ell},Y_{1},\ldots,Y_{\ell},K\},\\
[X_i,Y_i]=-b_i Z, \qquad~~~~~~~~~~~~~~~~~~~~~~ i=1,\ldots,\ell \\
[X_i,Y_j]=0,   \qquad~~~~~~~~ ~~~~~~~~~~~~~~~~~~~      i\neq j\\
[X_i,K]=[Y_i,K]=0,\qquad ~~~~~~~~~~~~~~~i=1,\ldots, \ell \\
[X_i,Z]=[Y_i,Z]=0,\qquad~~~~~~~~~~~~~~~~ i=1,\ldots, \ell 
\end{cases}
\eqn
where the new vector field $K$ is in the kernel of the bracket mapping, i.e. it commutes with all others elements.
Since abnormal extremals are never optimal in quasi contact case (see Remark \ref{r:abnormal}), we are reduced to compute the exponential map to find geodesics. With analogous computations of contact case we get the following expression for the exponential map from the origin \begin{align} \label{eq:exp4}
x_i(t)&=\dfrac{r_i}{b_iw}(\cos (b_iw t +\theta_i) - \cos \theta_i), \notag \\
y_{i}(t)&=\dfrac{r_i}{b_iw}(\sin (b_iw t +\theta_i) - \sin \theta_i), \\
x_{2\ell+1}(t)&=u_{2\ell+1} t, \notag\\
z(t)&=\dfrac{1}{2w^2}(|r|^2wt-\sum_i \dfrac{r_i^2}{b_i} \sin b_iw t). \notag
\end{align}

From \eqref{eq:exp4} it is easily seen that the jacobian of exponential map has exactly the same expression as in contact case \eqref{eq:exp3}. Since zero is always an eigenvalue of $L$, but is never the maximum one, we can proceed as in the contact case and all the regularity results extend to this case.

\section{Proof of Theorem \ref{t-nleq5} and extension to general corank 1 case}
We start this section with the proof of Theorem \ref{t-nleq5}, after that we extend the result to the general corank 1 case. 

\begin{proof}[Proof of Theorem \ref{t-nleq5}]
Let $\sub$ be a sub-Riemannian structure such that $\dim M\leq 5$.

\smallskip
$(i)$. If $\mc{G}(\sub)\neq (4,5)$.  From Theorem \ref{t:nf} we know that at every point $q\in M$, the nilpotent approximation $\widehat{\sub}_{q}$ has a unique normal form, hence by Corollary \ref{c:isometrie} all nilpotent approximations are isometric. From this property it easily follows that $f_{\mc{PS}}$, the Popp volume of the unit ball, is constant (recall that Popp measure is intrinsic for the sub-Riemannian structure). This also implies that for a smooth volume $\mu$ the density $f_{\mu \mc{S}}$ is smooth.

\smallskip
$(ii)$. If $\mc{G}(\sub)= (4,5)$ by Theorem \ref{t:nf} it is sufficient to consider the case when the family of nilpotent structure has the normal form \eqref{eq:45}, where $\alpha=\alpha(q)$ depends on the point. Notice that the formula \eqref{eq:drop} for the volume of the unit ball is still valid, where now $b_{1}=1$ and $b_{2}=|\alpha|$. 

Theorem \ref{t-c4c5} proves that the density is $\con^{3}$ at points where $|\alpha|>0$, i.e. in the contact case. We are then reduced to the study of the volume near a point where the eigenvalue $\alpha$ crosses zero. In particular we show that 
 the volume is smooth at these points. Since the eigenvalue $\al$ is approaching zero, it is not restrictive to assume $|\al|<1$. Let us consider then the function defined on the interval $(-1,1)$
\begin{equation} \label{eq:360}
W(\alpha)=\int_0^{\pi} \frac{1}{ \al^{2}s^{6}}\left(  \sin^2(\al s) \sin s(s\cos s-\sin s) +\sin^2 s \sin(\al s)(\al s\cos( \al s)-\sin(\al s))\right) ds.
\end{equation}
Note that $V(\alpha)=W(|\alpha|)$, where $V(\alpha)$ denotes the volume of the nilpotent ball relative to frequences $1$ and $\alpha<1$. It is easy to see that both
\bqn \label{eq:ultima}
\frac{\sin^{2} \alpha s}{\alpha^{2}} \qquad \text{and} \qquad  \frac{1}{\alpha^{2}}\sin(\al s)(\al s\cos( \al s)-\sin(\al s))
\eqn
are smooth as functions of $\alpha$ (also at $\alpha=0$). Hence $W$ is a smooth function (for $\al \in (-1,1)$).
Moreover it is easy to see that $W$ is an even smooth function of $\al$. Thus $W$ it is smooth also as a function of $|\alpha|$, which completes the proof.

\medskip
The same argument applies to prove that the $\con^{3}$ regularity holds in the general corank 1 case. Indeed, from \eqref{eq:explicit} and the fact that \eqref{eq:ultima} are smooth functions at $\al=0$, it follows that the integrand $G(q,s)$ is smooth as soon as one of the eigenvalue $b_{i}(q)$ is different from zero (recall that $b_{i}\geq0$ by definition). Since the structure is regular (i.e. the dimension of the flag do not depend on the point) and bracket generating, we have that $\max_{i} b_{i}(q)>0$ for every $q$, hence the conclusion. 
\end{proof}
\bigskip
\textbf{Acknowledgements.} We thank Gregoire Charlot who was the first to compute the volume of the nilpotent ball in the (4,5) case, and Jean-Paul Gauthier for very helpful discussions. We are also indebted with the anonymous reviewer for his valuable remarks.

{\small
\bibliography{bibliogr}
\bibliographystyle{abbrv}
}
\end{document}